\documentclass[final,5p,times,twocolumn]{elsarticle}

\usepackage{amsmath,amsfonts,amssymb,mathtools}
\usepackage{algorithmic}
\usepackage{algorithm}
\usepackage{array}
\usepackage[caption=false,font=normalsize,labelfont=sf,textfont=sf]{subfig}
\usepackage{textcomp}
\usepackage{stfloats}
\usepackage{url}
\usepackage{verbatim}
\usepackage{graphicx}
\usepackage{mathrsfs}
  
\usepackage{dsfont}
\usepackage{color}   
\usepackage{enumerate}
\usepackage{tabularx,arydshln}
\PassOptionsToPackage{hyphens}{url}\usepackage{hyperref}
\hypersetup{
    colorlinks = true,
    linkcolor = blue,
    urlcolor  = blue,
    citecolor = blue,
    anchorcolor = blue
}

\newcommand*{\QEDB}{\hfill\ensuremath{\square}}%
\DeclareMathOperator*{\esssup}{ess\,sup}

\newtheorem{theorem}{Theorem}
\newtheorem{proposition}{Proposition}
\newtheorem{lemma}{Lemma}
\newtheorem{corollary}{Corollary}
\newtheorem{definition}{Definition}
\newtheorem{assumption}{Assumption}
\newtheorem{remark}{Remark}
\newtheorem{example}{Example}
\newtheorem{problem}{Problem}

\journal{arXiv}

\begin{document}

\begin{frontmatter}

\title{\huge{Decentralized Concurrent Learning with Coordinated Momentum and Restart}}
\author[California]{Daniel E. Ochoa}
\author[Colorado]{Muhammad U. Javed}
\author[WashU]{Xudong Chen}
\author[California]{Jorge I. Poveda}

\address[California]{Department of 
Electrical and Computer Engineering, University of California, San Diego, La Jolla, CA 92093}	
\address[Colorado]{Department of 
Electrical, Computer and Energy Engineering, University of Colorado, Boulder, CO 80309}
\address[WashU]{Department of Electrical and Systems Engineering, Washington University in St. Louis, MO 63112}
\begin{abstract}
    This paper studies the stability and convergence properties of a class of multi-agent concurrent learning (CL) algorithms with momentum and restart. Such algorithms can be integrated as part of the estimation pipelines of data-enabled multi-agent control systems to enhance transient performance while maintaining stability guarantees. However, characterizing restarting policies that yield stable behaviors in decentralized CL systems, especially when the network topology of the communication graph is directed, has remained an open problem. In this paper, we provide an answer to this problem by synergistically leveraging tools from graph theory and hybrid dynamical systems theory. Specifically, we show that under a \emph{cooperative richness} condition on the overall multi-agent system's data, and by employing \emph{coordinated} periodic restart with a frequency that is tempered by the \emph{level of asymmetry} of the communication graph, the resulting decentralized dynamics exhibit robust asymptotic stability properties, characterized in terms of input-to-state stability bounds, and also achieve a desirable transient performance. To demonstrate the practical implications of the theoretical findings, three applications are also presented: cooperative parameter estimation over networks with private data sets, cooperative model-reference adaptive control, and cooperative data-enabled feedback optimization of nonlinear plants.
\end{abstract}

\end{frontmatter}

\section{INTRODUCTION}
Concurrent Learning (CL) techniques have emerged as powerful data-driven tools for designing estimation and learning dynamics in a wide range of applications where persistence of excitation (PE) conditions are either impractical or infeasible \cite{kamalapurkar2017concurrent,vamvoudakis2015asymptotically}. These techniques have demonstrated their utility in diverse fields, including parameter estimation in batteries \cite{ochoa2021accelerated}, exoskeleton robotic systems \cite{casas2022switched,casas2023switched}, mobile robots and aerial vehicles \cite{chowdhary2011theory}, extremum seeking algorithms \cite{poveda2021data}, and reinforcement learning controllers \cite{ochoa2022acceleratedADP,chowdhary2010concurrent}. In these applications, extensive datasets containing historical \emph{recorded} measurements of the relevant system signals are often available and can be leveraged for estimation purposes. When these datasets are ``sufficiently rich'', they can be seamlessly integrated into dynamic estimation algorithms, enabling (uniform) exponential convergence to the unknown parameters even in the absence of PE conditions.

However, relaxations of PE conditions can lead to suboptimal transient performance, particularly in terms of slow convergence rates that depend on the ``level of richness'' of the dataset used by the algorithm. Since datasets readily available in applications may exhibit prohibitively small levels of richness, there is a growing need for the development of enhanced CL techniques that can accelerate the convergence rate of the estimation dynamics while maintaining the desirable stability and robustness properties.

\subsection{Literature Review}
One promising direction to alleviate the slow convergence issue in decision-making algorithms is the incorporation of momentum with dynamic damping, see \cite{ochoa2021accelerated,le2022concurrent,nguyen2020momentumrnn}. For single-agent gradient-based dynamics with momentum, the use of decreasing damping has been shown to play a crucial role in inducing favorable acceleration properties \cite{muehlebach21, su2014differential, wibisono2016variational,nguyen2022improving}. However, it has also been shown that stability bounds in terms of $\mathcal{KL}$ functions may not exist for such systems unless the damping coefficients are persistently exciting \cite{poveda2020heavy}, a condition that precludes vanishing coefficients. Furthermore, it is well-known that, without proper tuning, the use of momentum can lead to undesirable oscillations \cite{o2015adaptive}. To address potential instability issues and to eliminate oscillatory behaviors, restart mechanisms that reset the momentum have been developed for single-agent systems using adaptive \cite{o2015adaptive,roulet2017sharpness} and periodic policies (usually called ``scheduled'') with carefully selected frequencies \cite{su2014differential,poveda2021robust,wang2022scheduled,o2015adaptive}. 
Recent works have also investigated the development of similar momentum-based algorithms in multi-agent systems, including distributed continuous-time heavy-ball dynamics with constant damping \cite{yu2020mass}, limiting equations of stochastic recursive algorithms as multi-agent flows with momentum \cite{boffi2020continuous}, and decision-making algorithms with momentum for high-order multi-agent systems \cite{sun2020continuous, ochoa2021momentum, ochoa2020robust}. However, existing approaches have primarily focused on \emph{undirected} network topologies. Additionally, the incorporation of momentum and restarting mechanisms in \emph{decentralized concurrent learning algorithms} has remained unexplored. Such algorithms are essential when a network of agents seeks to collaboratively and efficiently learn a common model by sharing local estimates with neighboring agents, without revealing their private data. Applications of these algorithms span various domains, including source seeking in autonomous mobile robots \cite{khong2014multi}, adaptive formation control of robotic teams \cite{chen2006smooth}, and cooperative adaptive control \cite{chen2013distributed}. 

\subsection{Contributions}
Motivated by the previous background, in this paper we study the synthesis and analysis of decentralized concurrent learning dynamics with momentum and restart for general directed graphs. In particular, we consider a model that extends the centralized dynamics studied in \cite{su2014differential,wibisono2016variational}, and \cite{poveda2021robust} to cases where each agent implements its own restart policy and shares information only with neighbors characterized by the topology of the communication graph. To assess the impact of the topology of graph on the stability properties of the dynamics, we exploit analytical tools from graph theory and hybrid dynamical system's theory \cite{goebel2012hybrid}. Using these tools, this paper makes the following primary contributions:

\textbf{(1)} We first introduce a class of multi-agent concurrent learning (CL) algorithms that incorporate momentum and a centralized restarting mechanism. We demonstrate that if: (a) the graph is strongly connected, (b) the overall data collected by the agents satisfies a ``cooperative richness condition," and (c) the restarting frequency exceeds a certain threshold that encodes the ``asymmetry" of the communication graph, then the resulting error estimation dynamics are input-to-state stable \cite{sontag1995characterizations} with respect to measurement noise and model error approximations. Furthermore, the convergence is exponential with rates assignable via the restarting period. These results are presented in Theorem \ref{theorem:centralized}.

\textbf{(2)} Next, by leveraging the robustness properties of the dynamics, we interconnect the momentum-based concurrent learning algorithms with a decentralized coordinated restarting mechanism, enabling a fully decentralized implementation. The resulting dynamical systems are also globally stable and exhibit convergence rates consistent with Theorem \ref{theorem:centralized} following an initial synchronization phase of the restarting times. These results are presented in Theorem \ref{theorem:decentralized}.

\textbf{(3)} Finally, we present three applications of the proposed momentum-based CL algorithms with restart within the context of data-enabled control: (a) cooperative parameter estimation without persistently exciting regressors in networks where nodes have private data with heterogeneous informativity properties; (b) data-enabled cooperative model-reference adaptive control; (c) data-enabled cooperative feedback-optimization. By employing (hybrid) Lyapunov-based techniques, we show that the resulting closed-loop systems exhibit favorable stability and convergence properties, which are also illustrated via numerical examples.

The rest of this paper is organized as follows: Section \ref{secpreliminaries} presents the preliminaries. Section \ref{secproblem} presents the problem formulation. Section \ref{section:centralized} presents the main results. Section \ref{section:applications} presents applications, Section \ref{sec:proofs} includes the proofs, and Section \ref{sec_conclusions} concludes the paper.
\section{Preliminaries}
\label{secpreliminaries}
\emph{Notation:} We use $r\mathbb{B}$ to denote a closed ball of appropriate dimension in the Euclidean space, of radius $r > 0$, and centered at the origin. Let $E_{ij}$ be the matrix with all entries equal to zero except the $ij^{\text{th}}$ entry, which is equal to one. Let  $\mathbf{1}_n\in\mathbb{R}^n$ be the vector of all ones, and $I_n\in\mathbb{R}^{n\times n}$ be the identity matrix. Given $x,y\in\mathbb{R}^n$, we let $(x,y)\coloneqq [x^\top,y^\top]^\top$ denote their concatenation. We use $\{e_1,e_2,\ldots,e_n\}$ to denote the standard basis of $\mathbb{R}^n$. A matrix $M\in \mathbb{R}^{N\times N}$ is represented in terms of its entries as $M=[m_{ij}]$, with $m_{ij}\in \mathbb{R}$ being its $ij^{\text{th}}$ entry. Similarly, we use $\mathbf{M} = [\mathbf{M}_{ij}]$ to represent a block matrix $\mathbf{M}$ in terms of its blocks, and use $\text{diag}\left(\{M_1,\ldots, M_J\}\right)$ to build a block diagonal matrix from the set of matrices $\{M_j\}_{j=1}^J$. Given a vector $x\in \mathbb{R}^n$ we let $\text{diag}(x)$ represent a diagonal matrix with diagonal entries $(i,i)$ given by the $i^{th}$ entry of $x$. A matrix $B=[b_{ij}]\in \mathbb{R}^{N\times N},$ is said to be nonnegative ($B\ge 0$) if $b_{ij}\ge 0$ for all $i,j$. The spectral radius of a matrix $B$ is denoted by $\rho(B)$. 
We use $|\cdot|$ for the vector 2-norm, $\|\cdot\|$ for the matrix norm induced by $|\cdot|$, and $|z|_{\mathcal{A}}\coloneqq \inf_{s\in\mathcal{A}}|z-s|$ to denote the distance of a vector $z\in\mathbb{R}^n$ to a closed set $\mathcal{A}\subset\mathbb{R}^n$. With a slight abuse of notation, given a matrix $\mathbf{P}\succeq 0$ we define the semi-norm $|x|_{\mathbf{P}}\coloneqq  (x^\top \mathbf{P}x)^{1/2}$. Given a set-valued mapping $M:\mathbb{R}^m\rightrightarrows \mathbb{R}^n$, the domain of $M$ is the set $\text{dom}(M)= \{x\in\mathbb{R}^m~:~M(x)\neq \emptyset\}$. A function $\gamma:\mathbb{R}_{\ge 0} \to\mathbb{R}_{\ge0}$ is of class $\mathcal{K}$ ($\gamma\in\mathcal{K}$) if it is continuous, strictly increasing, and satisfies $\gamma(0) =0$. It is said to be of class $\mathcal{K}_{\infty}$  ($\gamma\in\mathcal{K}_\infty$), if additionally $\gamma(r)\to\infty$ as $r\to\infty$. A function $\beta:\mathbb{R}_{\geq0}\times\mathbb{R}_{\geq0}\to\mathbb{R}_{\geq0}$ is of class $\mathcal{K}\mathcal{L}$ ($\beta\in\mathcal{KL}$) if it is nondecreasing in its first argument, nonincreasing in its second argument, $\lim_{r\to0^+}\beta(r,s)=0$ for each $s\in\mathbb{R}_{\geq0}$, and  $\lim_{s\to\infty}\beta(r,s)=0$ for each $r\in\mathbb{R}_{\geq0}$.

\emph{Graph Theory:} For a directed graph (digraph) $\mathcal{G} = (\mathcal{V},\mathcal{E})$, we denote by  $(i,j)\in \mathcal{E}$ a directed edge from node $i$ to node $j$, we call node $i$ an \emph{in-neighbor} of node $j$, and we call node $j$ an \emph{out-neighbor} of node $i$. We consider digraphs that do not have self-arcs. A weighted {\em Laplacian} matrix $\mathcal{L} = [l_{ij}] \in \mathbb{R}^{N\times N}$ associated with $\mathcal{G}$ satisfies the following: the off-diagonal entries are such that $l_{ij} < 0$ if $(i,j)$ is an edge, and $l_{ij} = 0$ otherwise; the diagonal entries $l_{ii}$ are determined such that every row of $\mathcal{L}$ sums to zero, and all its nonzero eigenvalues have positive real part \cite[Lemma 6.5]{bullo2018lectures}. A digraph is \emph{strongly connected} if for any two distinct nodes $i$ and $j$, there is a path from $i$ to $j$. The Laplacian matrix of a strongly connected digraph satisfies $\text{rank}(\mathcal{L}) = N-1$ \cite[Ch. 6]{bullo2018lectures}.

\emph{Hybrid Dynamical Systems:} In this paper, we work with dynamical systems that combine continuous-time and discrete-time dynamics. Such systems are called \emph{hybrid dynamical systems} (HDS) \cite{goebel2012hybrid}. The dynamics of a HDS with state $x\in\mathbb{R}^n$ are represented as:
\begin{subequations}\label{eq:HDS0}
\begin{align}
    &(x,u)\in C\coloneqq C_x\times\mathbb{R}^m,~~~~~\dot{x}\in F(x,u),\label{HDS0:flow}\\
    &(x,u)\in D\coloneqq D_x\times\mathbb{R}^m,~~~~x^+\in G(x), \label{HDS0:jump}
\end{align}
\end{subequations}
where $u\in\mathbb{R}^m$ is an exogenous input, $F:\mathbb{R}^n\times\mathbb{R}^m\rightrightarrows\mathbb{R}^n$ is called the flow map, $G:\mathbb{R}^n\rightrightarrows\mathbb{R}^n$ is called the jump map, $C\subset\mathbb{R}^n\times\mathbb{R}^m$ is called the flow set, and $D\subset\mathbb{R}^n\times\mathbb{R}^m$ is called the jump set. We use $\mathcal{H}=(C,F,D,G,u)$ to denote the \emph{data} of the HDS $\mathcal{H}$. These systems generalize purely continuous-time systems ($\tilde{D}=\emptyset$) and purely discrete-time systems ($\tilde{C}=\emptyset$). Time-varying systems can also be represented as \eqref{eq:HDS0} by using an auxiliary state $s\in\mathbb{R}$ with dynamics $\dot{s}>0$ and $s^+=s$. 
Solutions to system \eqref{eq:HDS0} are parameterized by a continuous-time index $t\in\mathbb{R}_{\geq0}$, which increases continuously during flows, and a discrete-time index $j\in\mathbb{Z}_{\geq0}$, which increases by one during jumps. Therefore, solutions to \eqref{eq:HDS0} are defined on \emph{hybrid time domains} (HTDs). Solutions to \eqref{eq:HDS0} are required to satisfy $\text{dom}(x)=\text{dom}(u)$, with $u(\cdot,j)$ being locally essentially bounded and Lebesgue measurable for each $j$. To establish this correspondence, a hybrid input $u$ in \eqref{eq:HDS0} is obtained from a suitable continuous-time input $u$ by using (with some abuse of notation) $u(t,j)=u(t)$  during the flows \eqref{HDS0:flow} for each fixed $j$, and by keeping $u$ constant during the jumps \eqref{HDS0:jump}. For a precise definition of \emph{hybrid time domains} and \emph{solutions} to HDS of the form \eqref{eq:HDS0}, we refer the reader to \cite{cai2009characterizations}.
To simplify notation, in this paper we use $|u|_{(t,j)}=\esssup_{\substack{(0,0)\preceq(\tilde{t},\tilde{j})\preceq (t,j)}}\left|u(\tilde{t},\tilde{j})\right|$, and we let $|u|_{\infty}\coloneqq\lim_{t+j\to\infty}|u|_{(t,j)}$.
%
%
The stability properties of HDS will be studied using the following notion.
\begin{definition}
Given a closed set $\mathcal{A}\subset C_x\cup D_x$, a HDS $\mathcal{H}$ of other form \eqref{eq:HDS0} is said to be \emph{input-to-state stable} (ISS) with respect to $|\cdot|_{\mathcal{A}}$ if there exist $\beta\in \mathcal{KL}$ and $\gamma\in \mathcal{K}$ such that every maximal solution pair $(x,u$) to $\mathcal{H}$ satisfies:  
\begin{equation}\label{KLbound}
|x(t,j)|_{\mathcal{A}}\leq \beta(|x(0,0)|_{\mathcal{A}},t+j) + \gamma(|u|_{(t,j)}),
\end{equation}
for all $(t,j)\in\text{dom}(x)$ and all $x(0,0)\in\mathbb{R}^n$. If \eqref{KLbound} holds with $u\equiv 0$, the set $\mathcal{A}$ is said to be \emph{uniformly globally asymptotically stable} (UGAS). If additionally, $\beta(r,s)=c_1re^{-c_2s}$ for some constants $c_1,c_2>0$, the set $\mathcal{A}$ is said to be \emph{uniformly globally exponentially stable} (UGES).  \QEDB
\end{definition}
%


\section{Problem Formulation}
\label{secproblem}
We consider a decentralized learning problem in a multi-agent system (MAS), where a group of $N\in\mathbb{Z}_{\ge 2}$ agents seeks to collaboratively estimate a common model characterized by a parameter $\theta^\star\in\mathbb{R}^n$. The agents share information with each other via a directed communication network modeled by a \emph{strongly connected} digraph $\mathcal{G}=\{\mathcal{V},\mathcal{E}\}$, where $\mathcal{V}\coloneqq\{1,2,\ldots,N\}$ is the set of nodes, and $\mathcal{E}$ is the set of edges. We assume that each agent $i\in\mathcal{V}$ has access to both real-time and past recorded measurements of a signal of the form
\begin{align}\label{outputsignal}
\psi_{i}^\star(t, d_i(t))=\phi_i(t)^\top \theta^\star + d_i(t),
\end{align}
where $d_i\in \mathbb{R}$ represents an unknown and possibly time-varying disturbance, and $\phi_i:\mathbb{R}_{\geq0}\to\mathbb{R}^n$ represents a regressor function (or basis functions), which is assumed to be continuous, uniformly bounded, and known to the $i^{th}$ agent. These assumptions are typical in single-agent  \cite{chowdhary2010concurrent,le2022concurrent,ochoa2021accelerated,vamvoudakis2015asymptotically} and distributed CL problems \cite{poveda2021data,chen2013distributed}. 
%
\subsection{Model Description and Key Assumptions}
To estimate $\theta^{\star}$, 
we consider a \emph{decentralized momentum-based concurrent learning} (DMCL) algorithm with the following update rule for the local estimate $\theta_i\in\mathbb{R}^n$ of each agent:
\begin{align*}
        \dot{\theta}_{i}(t)&=\frac{2}{\tau_i(t)}(p_i(t)-\theta_{i}(t)),~~~\dot{\tau}_{i}(t)\in [0,\omega],~~~\forall~~i\in\mathcal{V},
    \end{align*}
    where $\tau_i$ is a dynamic, non-decreasing coefficient, with rate of growth bounded by $\omega>0$, and which satisfies 
    $$\tau_i(t)\in[T_0,T],~~~\forall~t\in\mathbb{R}_{\geq0},~~~T>T_0>0,$$
    where $(T,T_0)$ are tunable parameters. The auxiliary state
    $p_i\in\mathbb{R}^n$ captures the incorporation of momentum, and it satisfies
\begin{align*}
        \dot{p}_i(t)=-2\tau_i(t) \left(\Lambda_i\left(\theta_i(t),\nu_i(t),t,\upsilon_i\right)+k_c\sum_{j\in\mathcal{V}}a_{ji}\left(\theta_i(t) - \theta_j(t)\right)\right),
    \end{align*}
where $k_c>0$ is a tunable gain, $\Lambda_i$ is a suitable mapping described below, and $a_{ji}$ is the $ji^{th}$ entry of the adjacency matrix of the graph $\mathcal{G}$ modeling the flow of information betwen agent $i$ and its neighbors.  The key components of the DMCL dynamics are explained below:
\begin{enumerate}[(a)]
\item The function $\Lambda_i$ has the general form:
\begin{equation}\label{maindynamics:individual:consensus}
\Lambda_i(\theta_i,\nu_i,t,\upsilon_i)=k_t\Psi_i\left(\theta_i,t,\upsilon_i\right)+k_r\Phi_i(\theta_i,\nu_i),
\end{equation}
where $k_{r}>0$ and $k_{t}\geq0$ are tunable constants. 
\item The function $\Psi_i$ in \eqref{maindynamics:individual:consensus} is given by
\label{eq:learningMaps}
\begin{align}
    \Psi_i(\theta_i,t,\upsilon_i(t)):=  \phi_i(t)\left(\hat{\psi}_i(\theta_i,t) - \psi_i^\star(t,\upsilon_{i}(t))\right),\label{psimapping}
\end{align}

\vspace{-0.4cm}
and it incorporates the \emph{real-time} information available to the $i^{\text{th}}$ agent, 
where $\psi_{i}^\star$ is given by \eqref{outputsignal},  $\hat{\psi}_i(\theta_i,t)\coloneqq \phi_i(t)^\top\theta_i$, and $\upsilon_i(t):=d_i(t)$ is a time-varying disturbance.

\item The function $\Phi_i$ in \eqref{maindynamics:individual:consensus} is given by
\begin{align}\label{eq:learningMaps2}
\hspace{-0.3cm}\Phi_i(\theta_i, \nu_i) \!\coloneqq\! \sum_{k=1}^{\bar{k}_{i}}\phi_i(t_{i,k})\!\!\left(\hat{\psi}_i(\theta_i,t_{i,k}) {-} \psi_i^\star(t_{i,k},\nu_{i,k})\right),
\end{align}

\vspace{-0.4cm}
and it incorporates past recorded measurements of the signal $\psi_{i}^\star$ in \eqref{outputsignal} and the regressor $\phi_i$, obtained at a sequence of times $\{t_{i,k}\}_{k=1}^{\bar{k}_i}$, where $\bar{k}_i\in\mathbb{Z}_{\geq0}$, and where $\nu_{i,k}:=d_i(t_{i,k})$ and 
$\nu_i \coloneqq (\nu_{i,1}, \nu_{i,2}, \dots, \nu_{i,\overline{k}_i}) \in \mathbb{R}^{\overline{k}_i}$ models the persistent disturbances occurring during data collection. 

\vspace{0.1cm}
\item The last term in the dynamics of $p_i$ captures the exchange of information between agent $i$ and its neighbors. Note that, in general, we have that $a_{ij}\neq a_{ji}$.
\end{enumerate}
To study the DMCL dynamics, the data matrix associated to the $i^{th}$ agent is defined as:
\vspace{-0.1cm}
\begin{align}\label{dataagenti}
\Delta_i \coloneqq \sum_{k_i=1}^{\bar{k}_i}\phi(t_{k_i})\phi(t_{k_i})^\top \in\mathbb{R}^{n\times n}.
\end{align}

\vspace{-0.2cm}
Rather than assuming that every matrix $\Delta_i$ is positive definite, as in standard single-agent concurrent learning (CL) \cite{chowdhary2010concurrent}, we will assume a weaker ``cooperative" richness condition on the overall data of the network \cite[Def. 2]{poveda2019codes}.
%
%
\begin{assumption}\label{assumption:CSR} 
There exists a constant $\alpha>0$, such that
\begin{equation}\label{CSRineq}
\sum_{i=1}^N \Delta_i \succeq \alpha I_n.
\end{equation}
Moreover, the graph $\mathcal{G}$ is strongly connected. \QEDB 
\end{assumption}

If \eqref{CSRineq} holds, the data $\{\Delta_i\}_{i\in\mathcal{V}}$ is said to be \emph{cooperatively sufficiently rich} (CSR).

\begin{remark}
 Assumption \ref{assumption:CSR} allows for some agents to have uninformative data (e.g., $\phi_i(t_{i,k})=0$) provided other agent's data is sufficiently rich data to satisfy \eqref{CSRineq}, see also \cite{javed2021excitation}. This is an important relaxation for large-scale MAS where, unlike standard CL \cite{chowdhary2010concurrent}, it might be unreasonable to assume that \emph{every} node's data satisfies $\Delta_i\succ0$. Moreover, note that in the DMCL dynamics agents do not share their data with other agents in the system. In fact, only the local estimates $\theta_i$ are shared with the neighboring agents. This prevents the direct solution of the estimation problem using ``single-shot" techniques, and instead calls for recursive algorithms that preserve the privacy of the individual data. \QEDB 
\end{remark}

\subsection{Connections to Accelerated Gradient Flows}
The form of the DMCL dynamics is closely related to the accelerated gradient flows with momentum studied in \cite{su2014differential,wilson2021lyapunov,wibisono2016variational,ochoa2021momentum}, which have the general form
\begin{align*}
\dot{x}_1(t)&=\frac{2}{\tau_c(t)}(x_2(t)-x_1(t)),\\
\dot{x}_2(t)&=-2\tau_c(t)\nabla f(x_1(t)),
\end{align*}
and where $f$ is a suitable convex cost function and $\tau_c:\mathbb{R}_{\geq0}\to\mathbb{R}_{>0}$ is a time-varying coefficient. Indeed, using the vectors $\theta:=(\theta_1,\theta_2,\ldots,\theta_N)$, $p:=(p_1,p_2,\ldots,p_N)$, the parameter error coordinates $\tilde{\theta}:=\theta-\mathbf{1}_N\otimes\theta^{\star}$, $\tilde{p}:=p-\mathbf{1}_N\otimes\theta^{\star}$, and the Laplacian matrix of the graph $\mathcal{L}$, the DMCL dynamics with a centralized coefficient $\tau=\tau_1=\ldots=\tau_N$ can be written as the following dynamical system:
\begin{equation}\label{errordynamicsmain}
\left(\begin{array}{c}
\dot{\tilde{\theta}}\\
\dot{\tilde{p}}
\end{array}\right)
=\hat{\mathbf{F}}(\tilde{\theta},\tilde{p},\tau,t),
\end{equation}
where $\hat{\mathbf{F}}$ is given by 
\begin{equation}\label{errormapping}
\hat{\mathbf{F}}(\tilde{\theta},\tilde{p},\tau,t)\!=\!\begin{pmatrix}
\dfrac{2}{\tau}(\tilde{p}-\tilde{\theta})\\
-2\tau \left(k_t\mathbf{A}(t)\!+\!k_r\mathbf{\Delta}\!+\!k_c\mathcal{\mathbf{L}}\right)\tilde{\theta}+\mathbf{U}(t)
\end{pmatrix}
\end{equation}
and where $\mathbf{L} \coloneqq \mathcal{L} \otimes I_n$, and where $\mathbf{A}$ and $\mathbf{\Delta}$ are block-diagonal matrices given by:

\vspace{-0.3cm}
\begin{small}
\begin{align*}
&\mathbf{A}(t)\coloneqq\text{diag}\left(\left\{\phi_1(t)\phi_1(t)^\top,\dots, \phi_N(t)\phi_N(t)^\top \right\}\right),\\
&\mathbf{\Delta}\coloneqq\text{diag}\left(\left\{\sum_{k=1}^{\bar{k}_1}\phi_1(t_k)\phi_1^\top(t_k),\dots, \sum_{k=1}^{\bar{k}_N}\phi_N(t_k)\phi_N^\top(t_k) \right\}\right)
\end{align*}
\end{small}
and $\mathbf{U}$ is given by:
\begin{small}
\begin{align}\label{defU}
\mathbf{U}(t)&:=\left[\begin{array}{c}
-2\tau k_t\phi_1(t)\upsilon_1(t)+k_c\sum_{k=1}^{\bar{k}_1}\phi_1(t_k)\nu_{1,k}\\
\vdots\\
-2\tau k_t\phi_{N}(t)\upsilon_N(t)+k_c\sum_{k=1}^{\bar{k}_N}\phi_N(t_k)\nu_{N,k}
\end{array}\right].
\end{align}
\end{small}
However, while similar decentralized algorithms have been studied in \cite{sun2020continuous, ochoa2021momentum, ochoa2020robust}, the DMCL dynamics do not describe a standard gradient flow with momentum due to the lack of symmetry on $\mathcal{L}$, i.e., the right-hand side of \eqref{errormapping} cannot be expressed as the gradient of a potential function, a property that usually plays a crucial role in the stability properties of momentum-based dynamics.  
\begin{figure*}[!ht]
\centering
\includegraphics[height = 3.8cm]{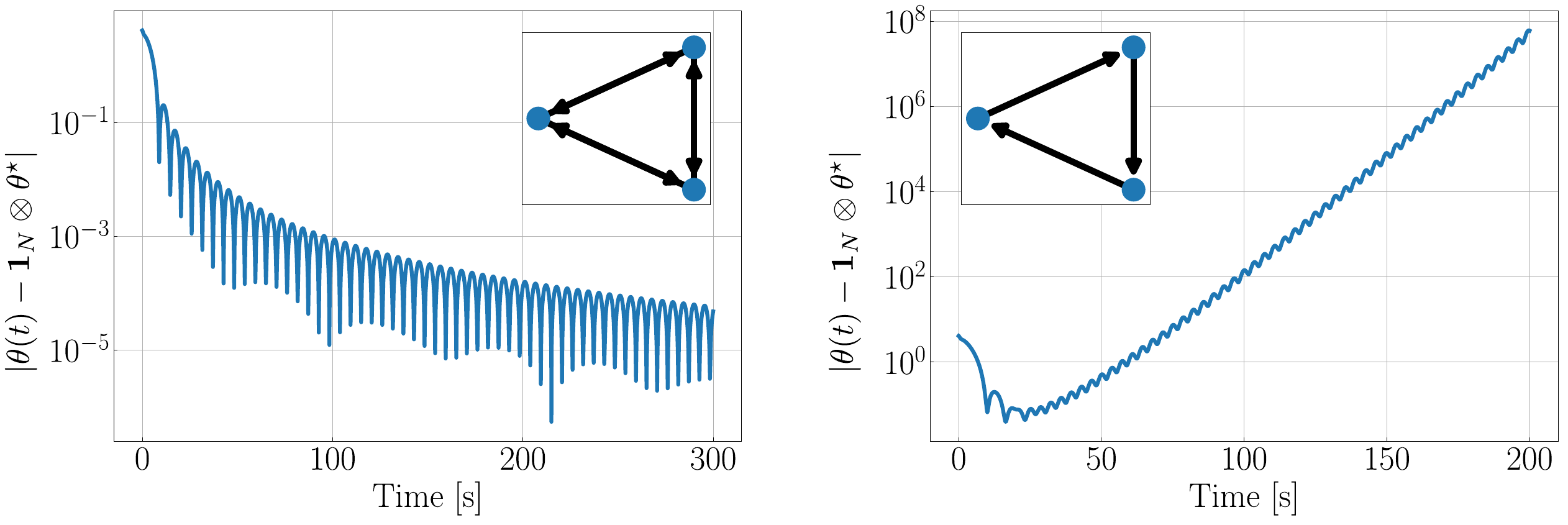}\hspace{0.02cm}
\includegraphics[height = 3.8cm]{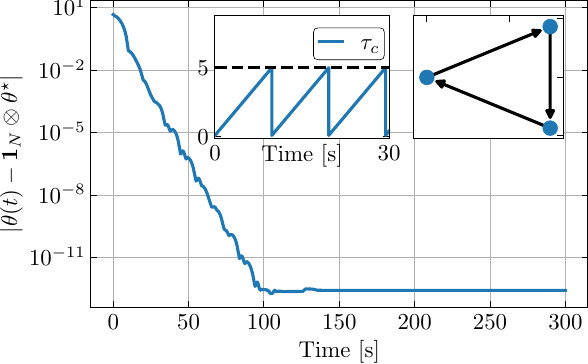}
\caption{Solutions to DMCL without restart can exhibit stability in symmetric graphs (left) and instability in asymmetric graphs (center). Stability in asymmetric graphs is recovered by employing a suitable coordinated restart mechanism (right). 
}
\label{fig:example}
\vspace{-0.2cm}
\end{figure*}
%
%
The following example highlights some of the challenges that can arise when momentum is used and the multi-agent system (MAS) has a communication topology characterized by a directed graph.

\vspace{0.1cm}
\begin{example}\label{example1}
Consider a multi-agent system with three agents, i.e., $\mathcal{V}=\{1,2,3\}$. We let $k_t=0$ and $d_i=0$, and for simplicity we assume that all agents use the same coefficient $\tau_c=\tau_1=\tau_2=\tau_3$, with $\tau(0)=0$ and $\omega=1/2$. We consider regressors $\phi_i(t)=(1,10e^{-it},100e^{-2it})$ with collected data satisfying Assumption \ref{assumption:CSR}, and the parameter $\theta^\star=(1,-2,1)$. The DMCL dynamics are implemented using $\dot{\tau}_c=1/2$. The left plot of Figure \ref{fig:example} shows the evolution in time (in logarithmic scale) of the estimation error $\tilde{\theta}=\theta-\mathbf{1}_N\otimes\theta^{\star}$ when the graph $\mathcal{G}$ is fully connected. As observed, the estimation error converges to zero, which is consistent with the stability results of \cite[Thm. 3]{su2014differential} and the fact that in this case, the DMCL dynamics describe an accelerated gradient system. Now, suppose that the communication graph is given by a directed cycle graph, as shown on the inset of the right plot of Figure \ref{fig:example}. In this case, the same DMCL algorithm ceases to be a momentum-based gradient flow and it exhibits the instability issue shown in the plot. The right plot, however, reveals a promising solution to the instability issue in asymmetric graphs. Stability can be restored by implementing a well-designed restart mechanism that accounts for the graph's asymmetry. The details of this mechanism will be elaborated upon in the following sections. \QEDB 
\end{example}

\vspace{-0.3cm}
\subsection{DMCL with Coordinated Restart}
To address the instability observed in Example \ref{example1}, while simultaneously inducing suitable convergence rates achieved via momentum, we can incorporate restart mechanisms into the algorithm. Such mechanisms persistently reset the momentum $\dot{\theta}_i$ and the dynamic coefficients $\tau_i$ whenever $\tau_i$ exceeds the upper bound $T$. The resets are performed according to the following discrete-time updates:
\begin{equation}\label{maindynamics:individual:jumps}
\theta_i^+=\theta_i,~~~p_{i}^+=p_{i}+\eta_i(\theta_{i}-p_{i}),~~~\tau_i^+=T_0,~~~\forall~i\in\mathcal{V},
\end{equation}
where $\eta_i\in\{0,1\}$ is a pre-defined parameter indicating the restart policy of agent $i$.
It has been shown that this approach can curtail the oscillations of momentum-based algorithms \cite{su2014differential,o2015adaptive} and also ``robustify'' their stability properties with respect to persistent disturbances \cite{poveda2021robust}. Indeed, note that the policy $\eta_i=1$ implies $p_{i}^+=\theta_i$, which implies $\dot{\theta}_i^+=0$. For multi-agent systems with undirected graphs, similar restart mechanisms of the form \eqref{maindynamics:individual:jumps} have been studied in \cite{ochoa2020robust,ochoa2021momentum}. However, the effectiveness of restarting in the context of multi-agent systems with \emph{directed} graphs has remained largely unexplored, and, as suggested by Example \ref{example1} the extension is non-trivial. Indeed, from the behavior observed in Figure \ref{fig:example} it should be clear that a ``slow'' restart policy (e.g., when $T$ is large) will not be able to stabilize the system. However, a ``fast'' restart policy (e.g., when $T$ is small) would keep $p_i$ approximately constant, thus reducing the effectiveness of using momentum. The right plot of Figure \ref{fig:example} shows the emerging behavior of the DMCL algorithm when restart is implemented by each node of the network with a ``suitable'' frequency and in a coordinated manner. While similar phenomena has been recently observed in game-theoretic problems \cite{ochoa2021momentum}, the use of momentum and restart in decentralized CL problems, and its dependence on the data of the system, the topology of the graph, and the perturbed models \eqref{outputsignal} have remained largely unexplored. This observation motivates the main research problem that we study in this paper:
%
%
\begin{problem}\label{problem1}
\normalfont{Characterize the restart mechanisms that: a) robustly stabilize the DMCL algorithm in directed networks; b) achieve ISS with respect to the disturbances $d_i$ in \eqref{outputsignal}; c) induce network-wide acceleration properties in the MAS.} \QEDB 
\end{problem}

\vspace{0.1cm}
In the next section we provide an answer to Problem \ref{problem1} using tools from hybrid dynamical systems and graph theory.
%
\section{Main Results}\label{section:centralized}
To tackle Problem 1, we first consider a centralized restart mechanism that makes use of a common state $\tau_c\in\mathbb{R}_{>0}$ that satisfies $\dot{\tau}_c\in[0,\omega]$. This centralized assumption simplifies the analysis and will be relaxed in the subsequent subsections to encompass decentralized implementations. For the purpose of analysis, we also use an auxiliary state $s\in\mathbb{R}$ with dynamics $\dot{s}=1$ to model any explicit dependence on time $t$.
\subsection{Centralized Restart: Hybrid Systems Model}
When using a common coefficient $\tau_c\in\mathbb{R}_{>0}$ to coordinate the restart of the DMCL algorithm, the resulting dynamical system can be modeled by the following differential inclusion, in vectorial form, with state $y_c\coloneqq(\theta,p,\tau_c,s)$:
%
%
\begin{align}\label{centralized:flowmap}
\dot{y}_c\!\in\!\mathbf{F}_c(y_c,u)\coloneqq
\begin{pmatrix}
\dfrac{2}{\tau_c}(p-\theta)\\
\!-2\tau_c\mathbf{\Lambda}(\theta, s,u)\!\\
[0,\omega]\\
1
\end{pmatrix},
\end{align}
whenever $y_c \in \mathbf{C}_c \times \mathbb{R}_{\ge 0}$, and where $u \coloneqq (\upsilon, \nu)$, the function $\mathbf{\Lambda}$ is given by
\begin{equation}\label{lamndamatrix00}
\mathbf{\Lambda}(\theta, s, u) \coloneqq k_t \mathbf{\Psi}(\theta, s, \upsilon) + k_r \mathbf{\Phi}(\theta, \nu) + k_c \mathbf{L} \theta,
\end{equation}
and the vectors $\upsilon$ and $\nu$ are defined as
\begin{align*}
\upsilon\coloneqq (\upsilon_{1}, \upsilon_{2}, \dots, \upsilon_{N}) \in \mathbb{R}^N,~~~\nu\coloneqq (\nu_{1}, \nu_{2}, \dots, \nu_{N}) \in \mathbb{R}^{\bar{k}},
\end{align*}
where $\bar{k}:=\sum_{i \in \mathcal{V}} \overline{k}_i$. The maps $\mathbf{\Psi}$ and $\mathbf{\Phi}$ above are defined as
\begin{subequations}\label{mapsdatatime}
\begin{align}
\!\!\!\!\!\!\!\mathbf{\Psi}(\theta, s, \upsilon)&\coloneqq (\Psi_{1}(\theta_1,s,\upsilon_{1}),\dots, \Psi_{N}(\theta_N,s,\upsilon_{N}))\\
\mathbf{\Phi}(\theta, \nu)&\coloneqq (\Phi_1(\theta_1,\nu_{1}),\dots, \Phi_N(\theta_N,\nu_{N})),\label{phidatamap} 
\end{align}
\end{subequations}
where the functions $\Phi_i,\Psi_{i}$ were defined in \eqref{psimapping}-\eqref{eq:learningMaps2} for all $i\in\mathcal{V}$. Since $\theta$ and $p$ are allowed to evolve in $\mathbb{R}^{nN}$, while $\tau_c\in[T_0,T]$, in \eqref{centralized:flowmap}, the flow set $\mathbf{C}_c$ is defined as:
\begin{equation}\label{centralized:flowset}
\mathbf{C}_c\coloneqq \mathbb{R}^{nN}\times \mathbb{R}^{nN} \times [T_0,T],
\end{equation}
which guarantees that during flows the state $\tau_c$ remains in the interval $[T_0,T]$. To incorporate restarts into the DMCL algorithm, each time $\tau_c$ meets the condition $\tau_c=T$, it is reset to $T_0$, and, all the states $(\theta_i,p_i)$ are updated as in \eqref{maindynamics:individual:jumps}. 
%
%
%
Therefore, using 
\begin{equation}
\mathbf{R}_{\eta}\coloneqq  \text{diag}(\eta)\otimes I_n,
\end{equation}
with $\eta=(\eta_1,\eta_2,\ldots,\eta_N)$, the discrete-time updates of the states $(\theta,p,\tau)$ of the hybrid system can be written in vectorial form as
\begin{align}\label{centralized:jumps}
\!\!\!(\theta^+,p^+,\tau_c^+)=\hat{\mathbf{G}}_c(x_c)=(\theta,p + \mathbf{R}_{\eta}(\theta {-} p),T_0) 
\end{align} 
which are executed whenever $(\theta,p,\tau)$ are in the jump set $\mathbf{D}_c$, defined as follows: 
\begin{equation}\label{jumpset}
\mathbf{D}_c\coloneqq  \mathbb{R}^{nN}\times \mathbb{R}^{nN} \times \{T\}.
\end{equation}
Therefore, the overall discrete-time dynamics of the system with state $y_c$ can be written as:
\begin{align}\label{discreteDMCL}
&y_c\in \mathbf{D}_c \times \mathbb{R}_{\ge 0},\quad 
y_c^+=\mathbf{G}_c(y_c)\coloneqq\hat{\mathbf{G}}_c(y_c)\!\times\!\{s\}.
\end{align}
By combining \eqref{centralized:flowmap} and \eqref{discreteDMCL}, the DMCL algorithm with centralized restart can be viewed as a HDS of the form \eqref{eq:HDS0}, with data
\begin{equation}\label{centralized:HDS}
\mathcal{H}_c\coloneqq  (\mathbf{C}_c\times \mathbb{R}_{\ge0},\mathbf{F}_c,\mathbf{D}_c\times \mathbb{R}_{\ge0},\mathbf{G}_c).
\end{equation}
Note that in this centralized HDS the jump set \eqref{jumpset} only imposes conditions on the state $\tau_c$. Namely, a restart is triggered whenever $\tau_c=T$. If $\dot{\tau}(t)=\text{constant}\in[0,\omega]$ for all time, then the HDS would model a DMCL algorithm with \emph{scheduled} periodic restart, where the time between two consecutive restarts is $(T-T_0)\omega^{-1}$. However, the differential inclusion in \eqref{centralized:flowmap} also allows us to consider scenarios where $\dot{\tau}$ is not constant but rather is any absolutely continuous function (between restarts) satisfying $\dot{\tau}\in[0,\omega]$, which includes functions that remain constant for arbitrarily long periods of time.

\begin{figure}[t]
    \centering
    \includegraphics[width = 0.99\linewidth]{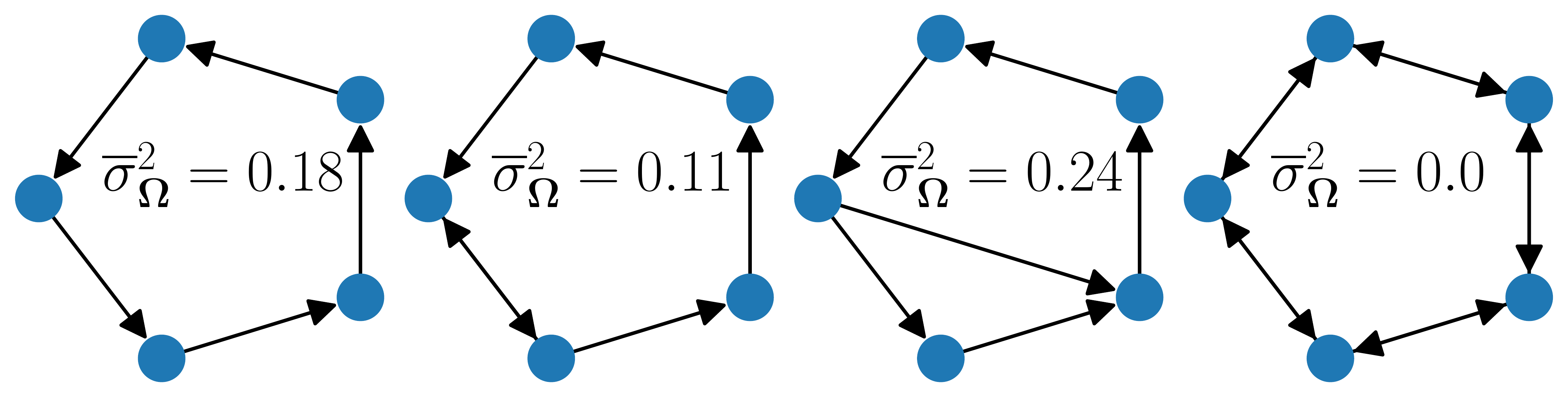}
    \caption{Parameter $\overline{\sigma}_{\mathbf{\Omega}}^2$ for strongly connected graphs with binary adjacency matrices and varying degrees of symmetry.}
    \label{fig:symmetryMeasureSigmaOmega}
\end{figure}

Before presenting our first main result, we introduce two technical propositions that play important roles in our results. All the proofs are presented in Section \ref{sec:proofs}.
\begin{proposition}\label{lemma:definitionQ}
Suppose that Assumption \ref{assumption:CSR} holds. Then, there exists a unit vector $q\in \mathbb{R}^N$ such that:
\begin{enumerate}[(a)]
\item The entries $q_i$ of $q$ satisfy:
\begin{align}\label{eq:Q_matrix}
    \overline{\sigma}_{\mathbf{Q} }\coloneqq  \max_{i\in\mathcal{V}} q_i \geq    \min_{i\in\mathcal{V}} q_i:=\underline{\sigma}_{\mathbf{Q} }>0.
\end{align}
\item $q^\top \mathcal{L} =0$ and $\mathcal{Q}\mathcal{L} + \mathcal{L}^\top \mathcal{Q}\succeq 0$ with $\mathcal{Q}\coloneqq \text{diag}(q)$.
\item The function $\mathbf{\Lambda}(\theta,s,u)$ in \eqref{lamndamatrix00} with $k_t=0$ and $ \nu=0$ can be decomposed as follows:
\begin{align}\label{decomposition1eq}
k_r \mathbf{\Phi}(\theta, 0) + k_c \mathbf{L} \theta= \mathbf{Q}^{-1}\left(\mathbf{\Sigma} + \mathbf{\Omega}\right)\tilde{\theta} ,
\end{align}
where $\mathbf{Q}\coloneqq \mathcal{Q}\otimes I_n$, $\tilde{\theta}:=\theta-\mathbf{1}_N\otimes\theta^{\star}$,
\begin{subequations}\label{eq:sigma_omega}
\begin{align}
\mathbf{\Sigma}&:=k_{r}\mathbf{Q}\mathbf{\Delta}\!+\! \frac{k_\text{c}}{2}\left(\textbf{Q}\textbf{L}\!+\!\textbf{L}^\top\textbf{Q}\right)\\
\mathbf{\Omega}&:=\! \frac{k_\text{c}}{2}\left(\textbf{Q}\textbf{L}\!-\!\textbf{L}^\top\textbf{Q}\right),
\end{align}
\end{subequations}
and $\mathbf{\Delta}\coloneqq \textbf{diag}\left(\{\Delta_1,\Delta_2,\dots, \Delta_N\}\right)$, where $\Delta_i$ is given by \eqref{dataagenti}. 
\item There exists a class-$\mathcal{K}_{\infty}$ function $\chi(\cdot)$ such that
\begin{align}\label{eq:omegabound}
\!\!\!\!\Bigl[ \mathbf{\Omega} + k_{t} \tilde{\mathbf{A}}(t)\biggr]   \biggl[ \mathbf{\Omega} + k_{t}   \tilde{\mathbf{A}}(t) \Bigr]^{\!\!\top} \!\!\!\preceq (\overline{\sigma}_{\mathbf{\Omega}}^2 + \chi(k_t)^2) I_{Nn},
\end{align}

\vspace{-0.2cm}
$\forall~t\geq0$, where $\tilde{\mathbf{A}}(t)\coloneqq\mathbf{Q}  \mathbf{A}(t)$ and $\overline{\sigma}_{\mathbf{\Omega}}$ is the largest singular value of $\mathbf{\Omega}$.  \hfill\QEDB
\end{enumerate}
\end{proposition}
\begin{remark}
By construction, if the Laplacian $\mathcal{L}$ is symmetric, then $\overline{\sigma}_{\mathbf{\Omega}}^{2} =0$. However, if $\mathcal{L}$ is asymmetric, then in general we have $\overline{\sigma}_{\mathbf{\Omega}}^{2} \neq 0$. For the purpose of illustration, Figure \ref{fig:symmetryMeasureSigmaOmega} presents four examples of different graphs $\mathcal{G}$ and their corresponding numerical values of $\overline{\sigma}_{\mathbf{\Omega}}^{2}$. \QEDB  
\end{remark}

\begin{proposition}\label{proposition2}
Suppose that Assumption \ref{assumption:CSR} holds; then, there exist $\overline{\sigma}_{\mathbf{\Sigma}}\geq \underline{\sigma}_{\mathbf{\Sigma}}>0$ such that 
\begin{equation}
\overline{\sigma}_{\mathbf{\Sigma}}I_{Nn}\succeq\mathbf{\Sigma}\succeq \underline{\sigma}_{\mathbf{\Sigma}} I_{Nn},
\end{equation}
where $\mathbf{\Sigma}$ is given by \eqref{eq:sigma_omega}. \QEDB 
\end{proposition}

\subsection{Input-to-State Stability of \texorpdfstring{$\mathcal{H}_c$}{Hc}}
%
%
%
%

%

\vspace{0.1cm}
With Propositions \ref{lemma:definitionQ}-\ref{proposition2} at hand, we are now ready to present the first main result of this paper, which provides conditions to stabilize the DMCL algorithm using a centralized restart. In particular, we study the stability properties of $\mathcal{H}_c$ with respect to the closed set $\mathcal{A}_c \coloneqq \mathcal{A}_{\theta p}\times [T_0,T] \times \mathbb{R}_{\ge 0}$, where
\begin{align}\label{centralized:stable_set}
\mathcal{A}_{\theta p}\coloneqq \{\mathbf{1}_N\otimes\theta^\star\}\times \{\mathbf{1}_N\otimes\theta^\star\},
\end{align}
which precisely describes the situation where all agent's estimates $\theta_i$ are equal to the true parameter $\theta^*$.

\vspace{0.1cm}
\begin{theorem}\label{theorem:centralized}
Suppose that Assumption \ref{assumption:CSR} holds, and let the constants ($\underline{\sigma}_{\mathbf{Q}}$, $\overline{\sigma}_{\mathbf{Q}}$, $\overline{\sigma}_{\mathbf{\Omega}}^2$, $\underline{\sigma}_{\mathbf{\Sigma}}$) be given by Proposition \ref{lemma:definitionQ}. If the restart parameters $(\omega,T_0,T)$ satisfy $\omega\in(0,1)$ and
\begin{align}\label{theorem:centralized:frequencyband}
\left(\frac{1}{2}\frac{\overline{\sigma}_{\mathbf{Q} }}{ \underline{\sigma}_{\mathbf{\Sigma}}}{+}{T_{0}^2}\right)^{\frac{1}{2}}\eqqcolon\underline{\mathbf{T}} < T< \overline{\mathbf{T}}\coloneqq\! \left(\frac{\underline{\sigma}_{\mathbf{Q} }(1-\!\omega)\underline{\sigma}_{\mathbf{\Sigma}} }{{ \overline{\sigma}^2_{\mathbf{\Omega}}+ \chi(k_t)^2   }}\right)^{\frac{1}{2}},
\end{align} 
%
%
then the following hold:
\begin{enumerate}[(a)]
\item For any restart policy $\eta\in\{0,1\}^N$ the HDS $\mathcal{H}_c$ renders the set $\mathcal{A}_c$ ISS with respect to the input $u$.

\item If $\eta_i=1$ for all $i\in\mathcal{V}$, and $\dot{\tau}_c\coloneqq \omega$, then, for every initial condition $y_0:=y_c(0,0)\in(\mathbf{C}_c\cup \mathbf{D}_c)\times\mathbb{R}_{\ge 0}$, every solution-input pair $(y_c,u)$ of $\mathcal{H}_c$, and every $(t_j,j)\in\text{dom}(y_c,u)$ with $t_j\coloneqq\min\{t:(t,j)\in\text{dom}(y_c,u)\}$, the sampled sequence of estimates $\theta(t_j,j)$ satisfies
\begin{equation}\label{theorem:centralized:convergencebound}
|\theta(t_j,j)-\mathbf{1}_N\otimes\theta^\star|^2\leq k_1\cdot\mu^j |y_0|_{\mathcal{A}_c}^2  + k_2|u|^2_{(t_j,j)},
\end{equation}
where $k_1,k_2>0$, and $\mu(T) \coloneqq (\underline{\mathbf{T}}/T)^2$. \hfill\QEDB
\end{enumerate}
\end{theorem}

The main result of Theorem \ref{theorem:centralized} 
reveals the impact of the asymmetry of $\mathcal{G}$ on the resetting parameter $T$. In particular, the following observations are in order:

\vspace{0.1cm}
\noindent
(1) When $\mathcal{L}$ is symmetric (i.e., $\overline{\sigma}_{\mathbf{\Omega}}^{2} =0$) and the DMCL dynamics do not use real-time data (i.e., $k_t=0$), condition \eqref{theorem:centralized:frequencyband} reduces to $\underline{\mathbf{T}}<T<\infty$, which can always be satisfied using \emph{any} positive constant $T$, recovering the results of \cite[Thm. 2]{poveda2021robust} in the context of optimization. 

\vspace{0.1cm}
\noindent
(2) In general, the more ``informative'' is the collective data in the overall system (i.e., the larger is $\alpha$ in \eqref{CSRineq}), the larger the parameter $\underline{\sigma}_{\mathbf{\Sigma}}$ will be, thus providing more flexibility to increase the upper bound $T$.

\vspace{0.1cm}
\noindent
(3) The ISS result implies that the trajectories of the algorithm will converge to a neighborhood of the true parameter $\theta^\star$, where the size of the neighborhood shrinks as the disturbances $d_i$ shrink in \eqref{outputsignal}. When $d_i=0$, the result establishes asymptotic convergence to the true parameter.

\vspace{0.1cm}
\noindent
(4) Lastly, when $u\equiv 0$, the bound \eqref{theorem:centralized:convergencebound} characterizes the ``accelerated'' convergence properties of $\mathcal{H}_c$ towards the true model $\theta^\star$. Since the rate of convergence is $T$-dependent, the rate of convergence depends on $T$. In particular, following similar steps as in the centralized case \cite{poveda2021robust}, the ``optimal'' value of $T$ that minimizes the contraction coefficient $\mu(T)$ over a given window of time can be computed as $T^*= e\left(\frac{\overline{\sigma}_{\mathbf{Q} }}{ 2\underline{\sigma}_{\mathbf{\Sigma}}}+{T_{0}^2}\right)^{\frac{1}{2}}$. 

\vspace{0.1cm}
Next, the following corollary leverages the expression of $T^*$ to obtain convergence bounds that parallel those obtained for centralized single-agent systems \cite{poveda2021robust}.

\vspace{0.1cm}
\begin{corollary}\label{corollary:centralized:optimalrestarting}
Suppose that all the assumptions of Theorem \ref{theorem:centralized} hold  with $T=T^*$, $\dot{\tau}_c=\omega$ and $u\equiv 0$; then, \eqref{theorem:centralized:convergencebound} holds with $k_2=0$, and for each $\varepsilon>0$ we have $|\theta(t_j,j)-\mathbf{1}_N\otimes\theta^\star|^2\leq\varepsilon$ for all $t_j>t_j^*$, where
$
t_j^*\coloneqq \frac{1}{2\omega}\left(T^*-T_0\right)\log \left(\frac{1}{\varepsilon} \frac{\overline{c}}{\underline{c}}|y_0|_{\mathcal{A}}^2\right).
$
\QEDB
\end{corollary}

\vspace{0.1cm}
The bound from Corollary \ref{corollary:centralized:optimalrestarting} implies that, as $T_0{\to}0^+$, the convergence of $\theta_i$ towards $\theta_i^*$ is of order $\mathcal{O}\left({e^{-\sqrt{\underline{\sigma}_{\mathbf{\Sigma}}/\overline{\sigma}_{\mathbf{Q}}}}}\right)$, for all $i\in\mathcal{V}$. We complete this section with a corollary for the case $\eta=0$,  which guarantees the ISS properties of $\mathcal{H}_c$, but not convergence bounds of the form \eqref{theorem:centralized:convergencebound}. 

\vspace{0.1cm}
\begin{corollary}\label{corollary:centralized:stabilityNopReset}
Suppose that Assumption \ref{assumption:CSR} holds, $\eta_i=0$ for all $i\in\mathcal{V}$, $\omega\in(0,1)$, and $ T_0 <  T< \overline{\mathbf{T}}$, with $\overline{\mathbf{T}}$ as defined in \eqref{theorem:centralized:frequencyband}. Then, the HDS $\mathcal{H}_c$ renders the set $\mathcal{A}_c$ ISS.\QEDB
\end{corollary}
The resetting bounds of Theorem \ref{theorem:centralized} and Corollary \ref{corollary:centralized:stabilityNopReset} only provide sufficient conditions for ISS with exponential convergence rates. 
It remains an open question how to obtain tight bounds on $(T_0,T)$ that are also \emph{necessary} for stability. We do not further pursue these questions in this paper.
\subsection{Decentralized Restart: Synchronization}\label{section:decentralized}
%
Since a central coordinator might not exist in large-scale networks, in this section, we study decentralized restart strategies based on each agent $i\in\mathcal{V}$ implementing an individual dynamic coefficient $\tau_i$. To simplify our discussion, in this section we assume that $\eta_i=1$ and $\dot{\tau}_i=\omega$ for all $i\in\mathcal{V}$, and that $k_t=0$, which allows us to remove the state variable $s$ and its associated dynamics. However, all our results can be extended to the case when time-varying regressors are included. 

When each agent implements its own coefficient $\tau_i$, the continuous-time DMCL dynamics \eqref{centralized:flowmap} become
\begin{equation}\label{decentralized:flowmap}
x\in \mathbf{C},~~\dot{x}=\mathbf{F}(x,u)\coloneqq \begin{pmatrix}
  2\mathcal{T}^{-1}(p-\theta)\\
-2\mathcal{T}(k_{r}\mathbf{\Phi}(\theta,u) + k_{c}\mathbf{L}\theta)\\
\omega \mathbf{1}_N  
\end{pmatrix},
\end{equation}
where the main state is now $x=(\theta,p,\tau)\in\mathbb{R}^{Nn}\times\mathbb{R}^{Nn}\times\mathbb{R}^N$, $\mathcal{T}\coloneqq\text{diag}(\tau\otimes \mathbf{1}_n)$, $\tau=(\tau_1,\tau_2,\dots,\tau_N)$, $\mathbf{\Phi}$ is given by \eqref{phidatamap}, and the flow set is given by
\begin{equation}\label{decentralized:flowset}
\mathbf{C}\coloneqq \mathbb{R}^N\times\mathbb{R}^N\times[T_0,T]^N.
\end{equation}
In this case, restarts of the form \eqref{maindynamics:individual:jumps} with $\eta_i=1$  occur whenever at least one of the agents satisfies the condition $\tau_i=T$. This behavior can be modeled by the following jump set:
\begin{equation}\label{decentralized:jumpset}
\mathbf{D}=\left\{x\in \mathbf{C}:~\max_{i\in\mathcal{V}}\tau_i=T\right\}.
\end{equation}
However, note that this approach would lead to uncoordinated restarts of the individual dynamics of the agents across the system.
%
%
For example, for any time window $[T_0,T]$, one can select $N$ equidistant initial conditions $\tau_i(0,0)\in[T_0,T]$, where $i\in\{1,2,\ldots,N\}$, which result in solutions experiencing $N$ restarts during this time window, each restart separated by intervals of flow of length $\frac{T-T_0}{N}$. Therefore, as $N\to\infty$, asynchronous restarts would occur more often, hindering the advantages of incorporating momentum into the flows of the algorithm to accelerate the overall system. 

To address this issue, and inspired by the synchronization algorithms of \cite{javed2021scalable}, we integrate the restart dynamics \eqref{maindynamics:individual:jumps} of each agent with a decentralized coordination mechanism for the states $\tau_i$. Specifically, each agent $i\in\mathcal{V}$ performs individual restarts of the form \eqref{maindynamics:individual:jumps} when $\tau_i=T$. However, the agents also implement the following \emph{additional} discrete-time updates whenever their neighbors $j\in\mathcal{N}_i$ satisfy the condition $\tau_j=T$:
\begin{equation}\label{decentralized:coordinationmap}
\tau_i^+\in \mathcal{R}_i(\tau_i)\coloneqq \left\{\begin{array}{cl}
T_0&\text{if}~\tau_i\in[T_0,r_i)\\
\{T_0,T\}&\text{if}~~\tau_i=r_i\\
T&\text{if}~\tau_i\in(r_i,T]
\end{array}\right.,
\end{equation}
where $r_i>0$ is a tunable parameter. To incorporate these additional discrete-time updates into the overall jump map of the system, consider the set-valued map
{\small
\begin{align*}
G_{d}(x)\coloneqq \Big\{(\hat{\theta},\hat{p},\hat{\tau})\in&~\mathbb{R}^{(2n+1)N}:\hat{\theta}=\theta,~
\hat{p}_i=p_i,~\hat{\tau}_i=T_0,\\
&\tau_j\in \mathcal{R}_j(\tau_j),~\hat{p}_j=p_{j},~\forall j\in \mathcal{N}_i,\\
&\hat{p}_k=p_k,~\hat{\tau}_k=\tau_k~\forall k\neq i\neq j\Big\},
\end{align*}
}
which is defined to be non-empty if and only if $\tau_i=T$ and $\tau_j\in[0,T)$. In words, the mapping $G_{d}(x)$ captures the resets of the \emph{individual states} $(\theta_{i},p_{i},\tau_i)\in\mathbb{R}^{2n+1}$ of agent $i$ via \eqref{maindynamics:individual:jumps}, and also the updates of its neighbors $j\in\mathcal{N}_i$ via \eqref{decentralized:coordinationmap}. The overall jump-map of the multi-agent hybrid system can then be defined using the outer-semicontinuous hull of $G_{d}$\footnote{The outer-semicontinuous hull of a set-valued mapping $G:\mathbb{R}^n\rightrightarrows\mathbb{R}^n$ is the unique set-valued mapping $\bar{G}:\mathbb{R}^n\rightrightarrows\mathbb{R}^n$ satisfying $\text{graph}(\bar{G})=\text{cl}(\text{graph}(G))$, where $\text{cl}(\cdot)$ stands for the closure.}, denoted $\overline{G}_{d}$ leading to
\begin{equation}\label{discretetimedecentralized}
x\in \mathbf{D},~~~x^+\in \mathbf{G}(x)\coloneqq \overline{G}_d(x).
\end{equation}
Note that system \eqref{discretetimedecentralized} preserves the sparsity property of the graph $\mathcal{G}$. 

The decentralized continuous-time dynamics \eqref{decentralized:flowmap} and the decentralized discrete-time dynamics \eqref{discretetimedecentralized} comprise the overall DMCL algorithm with restarts studied in this paper. This algorithm is fully modeled by the HDS 
\begin{equation}\label{mainalgorithmDMCL}
\mathcal{H}\coloneqq (\mathbf{C},\mathbf{F},\mathbf{D},\mathbf{G})
\end{equation}
with state $x=(\theta,p,\tau)\in\mathbb{R}^{Nn}\times\mathbb{R}^{Nn}\times\mathbb{R}^N$.

The following theorem provides a decentralized version of Theorem \ref{theorem:centralized}. In this case, stability of $\tau$ is studied with respect to the ``synchronized" set $\mathcal{A}_{\text{sync}}\coloneqq [T_0,T]\cdot\mathbf{1}_N\cup \{T_0,T\}^N$, and the stability properties of the overall state are studied with respect to
\begin{equation}\label{mainsetstability}
\mathcal{A}\coloneqq \mathcal{A}_{\theta p}\times \mathcal{A}_{\text{sync}}.
\end{equation}
For simplicity, we state the result for the case $u=0$, but we comment on the robustness properties of the dynamics.
\begin{theorem}\label{theorem:decentralized}
Consider the HDS $\mathcal{H}$ and suppose that Assumption \ref{assumption:CSR} holds and that:
\begin{enumerate}[(a)]
\item The parameters $(T_0,T)$ satisfy \eqref{theorem:centralized:frequencyband}.
\item The constants $\{r_i\}_{i\in \mathcal{V}}$ satisfy $T_0<r_i<T_0+\frac{(T-T_0)}{N-1}$
\end{enumerate}
Then, the set $\mathcal{A}\coloneqq \mathcal{A}_{\theta p}\times \mathcal{A}_{\text{sync}}$ is UGES for $\mathcal{H}$, and there exists a time $t^*\in \left[0,2\frac{T-T_0}{\omega}\right)$ such that for every solution $x$ of $\mathcal{H}$ and every $(t,j)\in\text{dom}(y)$ such that $t+j\geq t^*+2N$, the bound \eqref{theorem:centralized:convergencebound} holds.\QEDB
\end{theorem}
\begin{figure*}[t!]
    \centering
    \includegraphics[width=0.45\linewidth]{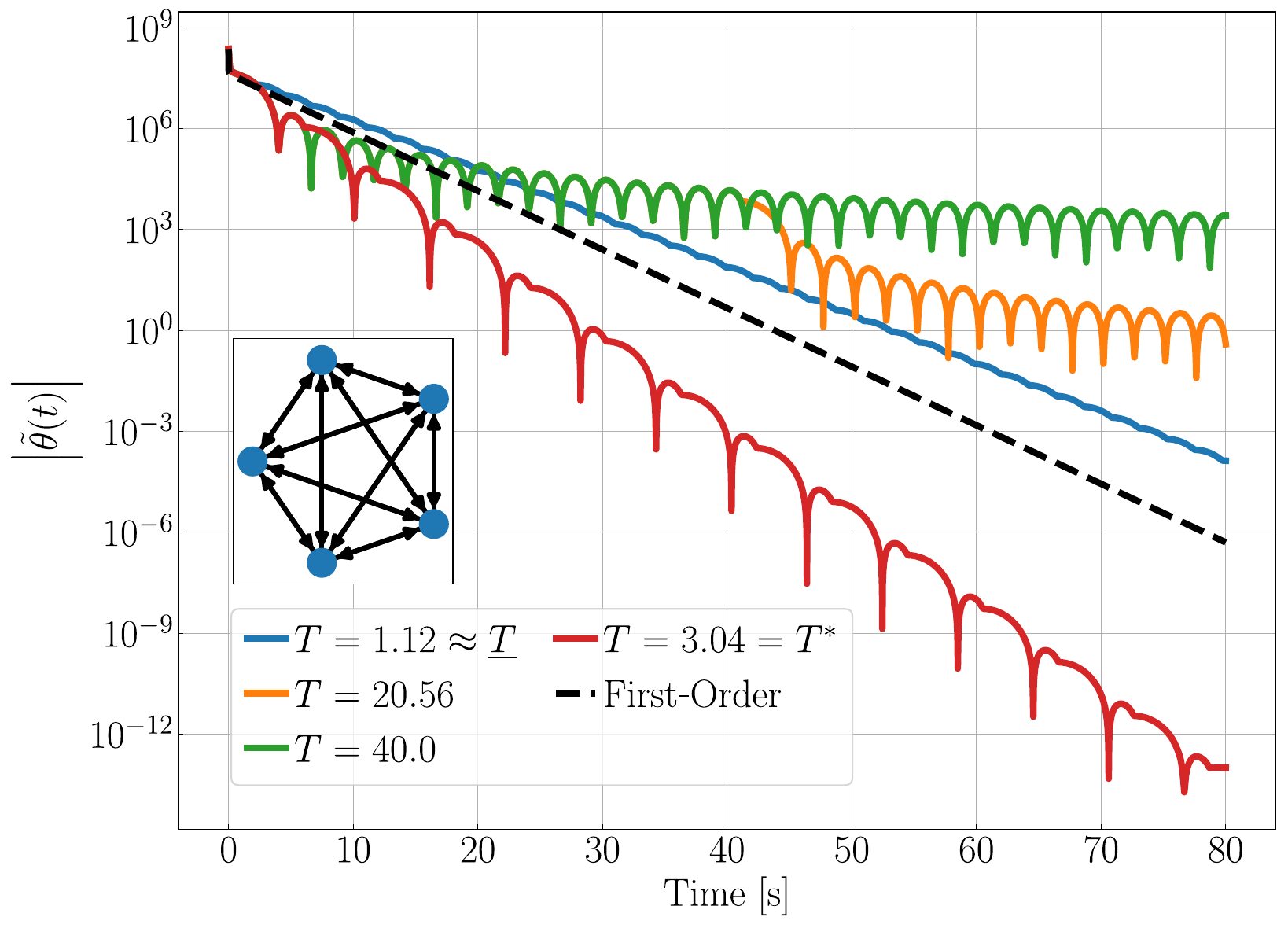}
    \includegraphics[width=0.45\linewidth]{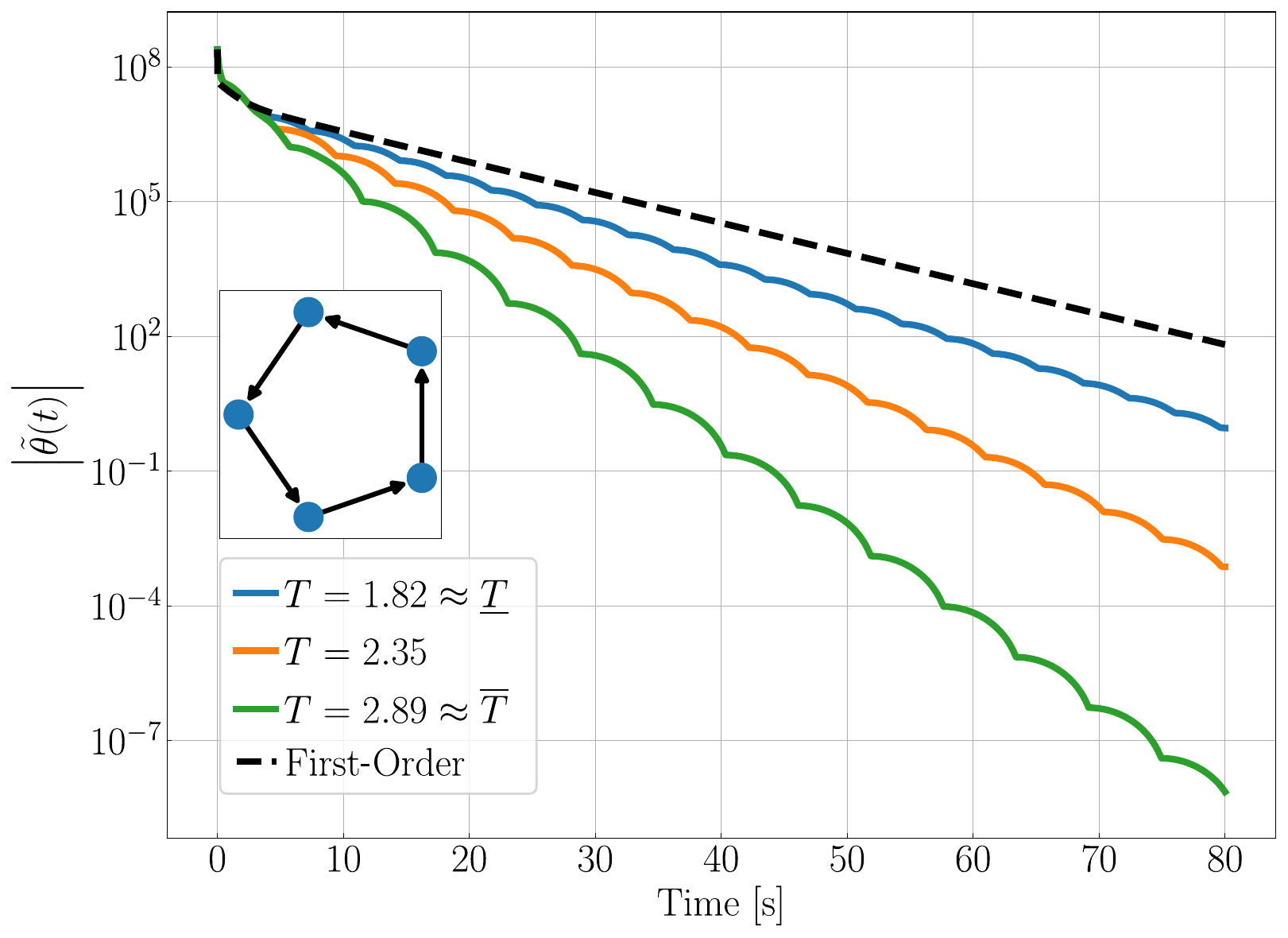}
    \caption{Left: Trajectories of $\mathcal{H}$ when $\mathcal{G}$ is fully connected. Right: Trajectories of $\mathcal{H}$ when $\mathcal{G}$ is a cycle. Here, $\tilde{\theta}=\theta-\mathbf{1}_N\otimes\theta^\star$}
    \label{fig:cooperativeIdentificationOverGraphs}
\end{figure*}
\begin{remark}[Nominal Robustness]\normalfont
 Since the hybrid system $\mathcal{H}$ is nominally well-posed in the sense of \cite[Ch. 6]{goebel2012hybrid}, the UGES properties of the DMCL algorithm are preserved, in a semi-global practical sense, under arbitrarily small additive perturbations on states and dynamics. This property is crucial for the use of $\mathcal{H}$ in practical applications where dynamic disturbances are unavoidable. \QEDB
\end{remark}
\begin{remark}[Strong Robustness via ISS]\normalfont
The techniques employed to proof Theorem \ref{theorem:decentralized} can be further utilized to obtain  ISS of $\mathcal{H}$ provided that $u$ originates from a dynamical system evolving in a compact set. We omit this extension due to space constraints.\QEDB
\end{remark}
\begin{remark}\normalfont
Since system $\mathcal{H}$ has no finite-escape times due to the global Lipschitz property of $\mathbf{F}$ in $\mathbf{C}$, it follows that the stability results of Corollary \ref{corollary:centralized:stabilityNopReset} also extend to $\mathcal{H}$ with $\eta_i=0$ for all $i$, recovering the convergence result of Theorem \ref{theorem:centralized} after an initial finite synchronization phase. \QEDB
\end{remark}

To the best of the author's knowledge, Theorems \ref{theorem:centralized}-\ref{theorem:decentralized} and the respective corollaries, are the first stability results for momentum-based CL algorithms implemented in multi-agent systems with general directed graphs. We note that in the literature of centralized CL, other accelerated algorithms have been studied using finite-time and fixed-time stability tools \cite{ochoa2021accelerated,rios2017time,tatari2021fixed}. However, as shown in the comparison presented in \cite{ochoa2021accelerated}, when the ``level of richness'' of the data (i.e., $\alpha$ in Assumption 1) is ``low'', momentum-based methods can yield better transient performance compared to other first-order non-smooth techniques. For \emph{decentralized} problems defined over networks, we are not aware of finite-time or fixed-time CL algorithms that are stable under Assumption \ref{CSRineq}. A natural progression for future research involves developing such algorithms and comparing them with the DMCL algorithms proposed in this paper.
\section{Applications in Estimation, Control, and Model-free Feedback Optimization}
\label{section:applications}
In this section, we apply the DMCL algorithm with restart in three different applications.
\subsection{Hybrid Cooperative Identification Over Digraphs}
First, we validate Theorem \ref{theorem:decentralized} in an cooperative estimation problem defined in a multi-agent system with $N=5$, $n=3$, and $\psi_{i}(s)=(10e^{-i s}-1)^2$, for all $i\in\mathcal{V}$. To implement the DMCL algorithm with coordinated restarts, we parameterize $\psi_{i}(\cdot)$ using the regressor $\phi_i(s)\coloneqq (1,10e^{-i s},100e^{-2is})$ and $\theta^\star=(1,-2,1)$. To satisfy Assumption \ref{assumption:CSR} with $\alpha = 5.5$, each agent records five measurements of $\psi_i$. We implement the hybrid system $\mathcal{H}$ and plot the resulting trajectories of the estimation error in the left plot of Figure \ref{fig:cooperativeIdentificationOverGraphs}, using $k_{r}=80, k_c=0.08$, and a fully connected graph. We also show with dashed lines the trajectory obtained when using the first-order decentralized CL dynamics of \cite{poveda2021data}. Since the graph is symmetric, in this case $T$ can be selected arbitrarily large to tune the convergence rate of the dynamics (see inequality \eqref{theorem:centralized:frequencyband}). The simulations start from a non-synchronized initial condition $\tau(0,0)\neq \tau_0 \mathbf{1}_5$ and rapidly achieve synchronization. Trajectories related to different choices of $T$ are also shown to illustrate the impact of the restart period on the convergence rate. Next, we let $\mathcal{G}$ be a cycle digraph, for which $\overline{\sigma}_{\mathbf{\Omega}}^2=0.18$. The resetting parameter $T$ is selected to satisfy inequality \eqref{theorem:centralized:frequencyband}, and the resulting trajectories are shown in the right plot of Figure \ref{fig:cooperativeIdentificationOverGraphs}. In this case, the best transient performance is obtained as $T$ approaches the upper bound $\overline{\mathbf{T}}$.
%
%
%
%
%
%
%
\subsection{Data-Enabled Hybrid Cooperative MRAC}
A key advantage of the robust stability results presented in Theorems 1 and 2, is that the DMCL dynamics can be interconnected with other systems for the solution of feedback control problems. To illustrate this application, we consider a multi-agent dynamical system, where each agent has individual dynamics of the form:
\begin{equation}\label{mrac}
\dot{\chi}_i=A_i\chi_i+B_iu_i + B_i\tilde{\psi}_i(\theta^\star,\chi_i),~~\chi_i\in\mathbb{R}^n,~u_i\in\mathbb{R}^m,
\end{equation}
where $\tilde{\psi}_i(\theta^\star,\chi)=\phi_i\left(\chi\right)^\top \theta^\star$ models structured uncertainty parameterized by a common vector $\theta^\star$, and a regressor $\phi_i$ that is known by each agent $i$. The agent's goal is to be able to asymptotically track a common bounded reference $r$ despite the uncertainty in their model.

\begin{figure*}[t]
\includegraphics[width=0.46\linewidth]{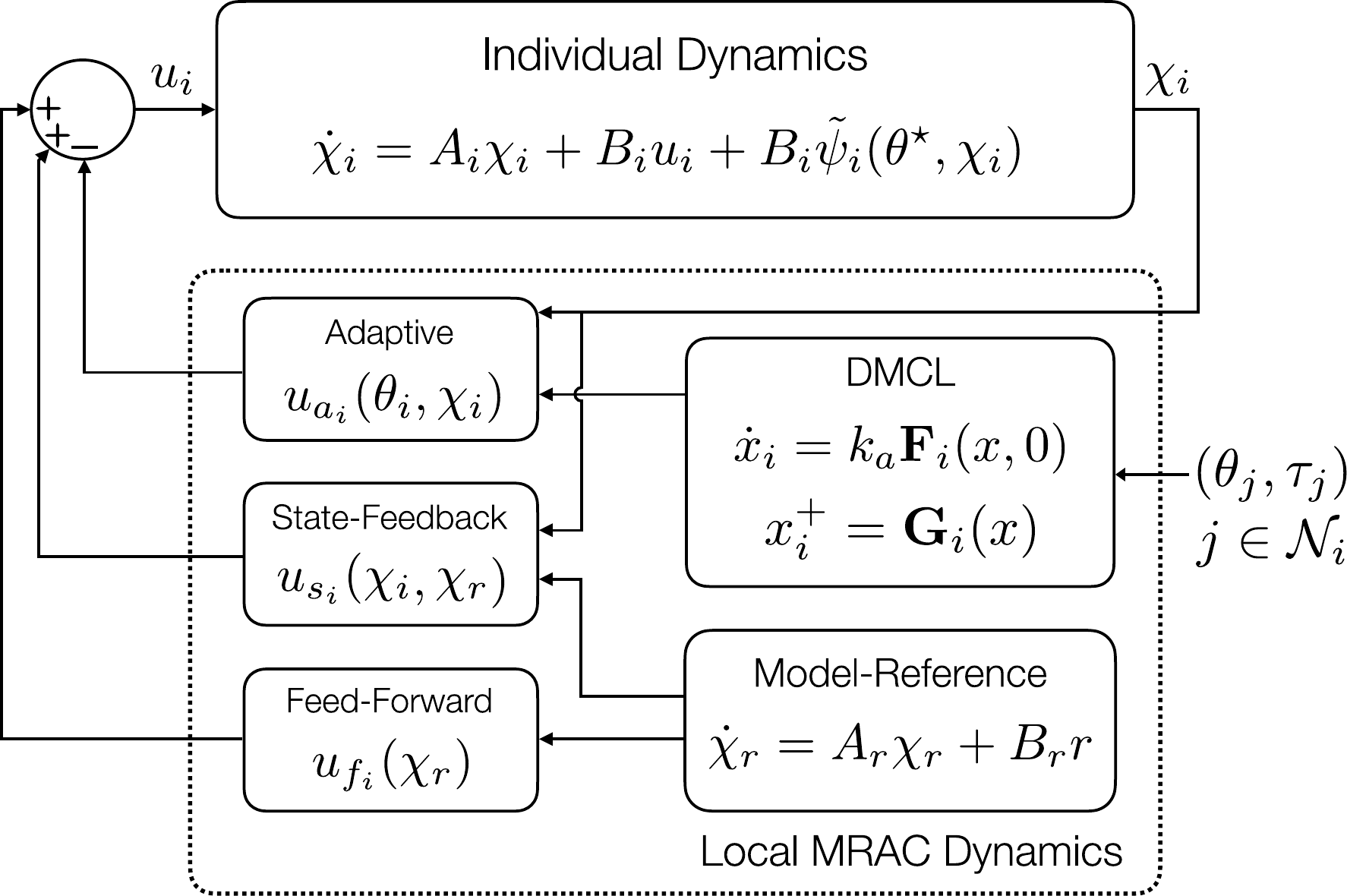}
\includegraphics[width=0.51\linewidth]{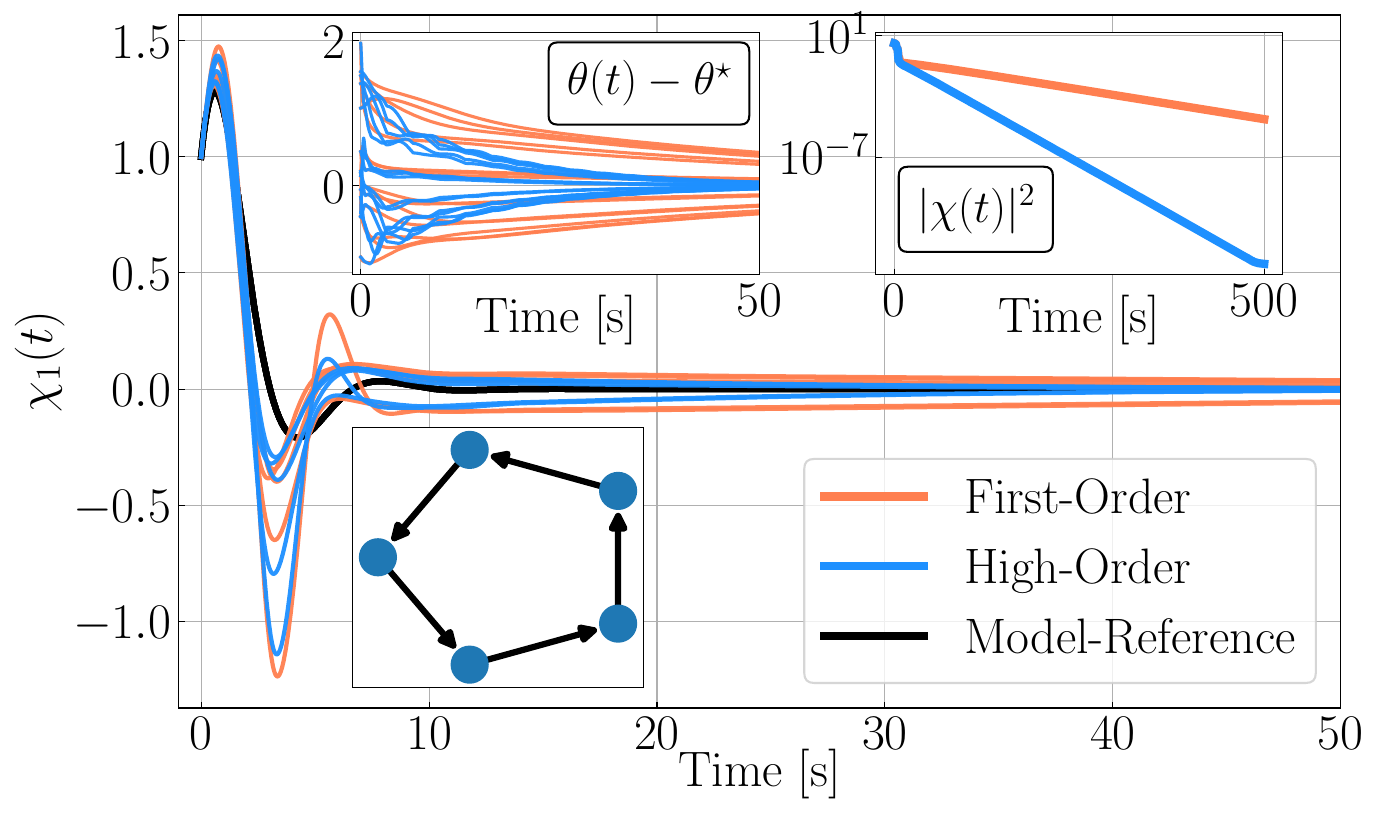}
    \caption{Left: Scheme of the $i^{th}$ agent's dynamics in the Cooperative MRAC. $\mathbf{F}_i$ and $\mathbf{G}_i$ represent the components of the overall DMCL flow-map and jump-map, respectively, corresponding to the state $x_i=(\theta_i,p_i,\tau_i)\in\mathbb{R}^{2n + 1}$. Right: Trajectories resulting from the Cooperative MRAC when $N=5$.}
    \label{fig:hybridMRAC}
\end{figure*}

\vspace{0.1cm}
\subsubsection{Two-Time Scale Hybrid Dynamics}
To solve the tracking problem we use a two-time scale approach. First, we introduce a reference model $\dot{\chi}_r=A_r \chi_r+B_rr$, where $A_r$ is assumed to be Hurwitz. Following the ideas of \cite{chowdhary2010concurrent}, each agent implements a model-reference adaptive control (MRAC) law that incorporates three elements: (1) an adaptive component $u_{a_i}(\theta_i,\chi_i)=\phi_i(\chi_i)^\top\theta_i$, where $\theta_i$ is the individual estimate of $\theta^\star$; (2) a state-feedback component $u_{s_i}(\chi_i,\chi_r)=-K(\chi_i-\chi_{r})$; and (3) a feed-forward term $u_{f_i}$ designed such that $B_iu_{f_i}(\chi_r)=(A_r-A_i)\chi_{r}+B_r r$; see Figure \ref{fig:hybridMRAC} for an illustration of the control law. Using $u_i(\theta_i,\chi_i,\chi_r) = u_{s_i}(\chi_i,\chi_r) + u_{f_i}(\chi_r) - u_{a_i}(\theta_i,\chi_i) $, and the error coordinates
%
%
%
$e_i=\chi_{i}-\chi_{r}$, the error dynamics for agent $i$ become:
%
\begin{equation}\label{errordynamics_agenti}
\dot{e}_i=A_{m_i}e_i+B_i\left(\tilde{\psi}_i(\theta^\star, e_i + \chi_r)-u_{a_i}(\theta_i,e_i+\chi_r)\right),
\end{equation}
where $A_{m_i}\coloneqq A_i - B_iK$, for all $i$. We make the assumption that that system \eqref{errordynamics_agenti} has no finite escape times from all initial conditions, and that Assumption \ref{assumption:CSR} holds. To cooperatively estimate $\theta$, we interconnect \eqref{errordynamics_agenti} with the DMCL algorith with restart given by \eqref{mainalgorithmDMCL} with flow map 
\begin{equation}\label{flomapnormalized}
x\in \mathbf{C},~~~~\dot{x}=k_a\mathbf{F}(x,0),    
\end{equation}
where the pair $(\mathbf{C},\mathbf{F})$ is given by \eqref{decentralized:flowmap}, and where $k_a>0$ is a tunable parameter.
%
%

To study the stability of the interconnected system, we will first analyze the system under a centralized timer $\tau_c$, and we interpret the closed-loop system as a two-time scale hybrid dynamical system with the DMCL algorithm having continuous-time dynamics operating in a faster time scale compared to \eqref{errordynamics_agenti}. Since $A_m$ is Hurwitz, for each $Q\succ0$ there exists $P\succ0$ such that $A_m^\top P +PA_m=-Q$, i.e., system  \eqref{errordynamics_agenti} is UGES when $\theta_i=\theta$. Similarly, by Theorem \ref{theorem:centralized}, the momentum-based hybrid dynamics $\mathcal{H}_c$ render the set $\mathcal{A}_c$ UGES via a Lyapunov function $V$. We can then study the interconnection of both systems using the Lyapunov function $V_1=0.5\tilde{V}(e)+ 0.5V(x)$, with $\tilde{V}(e)=e^\top P e$. From the proof of Theorem \ref{theorem:centralized}, during jumps $V_1$ satisfies $\Delta V\leq0$ because $e^+=e$. On the other hand, during flows of the closed-loop system the function $V$ satisfies
%
%
%
%
\begin{align}
\dot{V}&=-e^\top Q e - 0.5kV(y_c)+e^\top Q\phi(\chi(t))^\top \theta\\
&\leq -\lambda_{\min}(Q)|e|^2- k_a|y_c|_{\mathcal{A}_c}^2+k\overline{\phi}|e||y|_{\mathcal{A}},
\end{align}
where we used the quadratic lower bounds of $V$, and the boundedness of the regressors to obtain $k\overline{\phi}>0$. From here the result follows by completing squares and taking $k_a$ sufficiently large such that $\dot{V}<0$. Since $\Delta V\leq 0$, and the jumps are periodic, we obtain UGES of the set $\mathcal{A}_1=\{0\}\times\mathcal{A}$, where $\mathcal{A}$ is given by \eqref{mainsetstability}. The stability properties for the decentralized case follow now by leveraging the absence of finite escape times, and by using the reduction principle as in the proof of Theorem \ref{theorem:decentralized}.
%
%

\subsubsection{Numerical Example}
To illustrate the previous result, we consider a multi-agent system with $N=5$ agents, where the communication graph $\mathcal{G}$ is a directed cycle graph, see the inset in Figure \ref{fig:hybridMRAC}. We consider open-loop unstable individual dynamics characterized by matrices
$A_i= E_{12}\in\mathbb{R}^{2\times 2},\quad 
B_i = \left(
    0, \frac{2i - 1}{2i}
\right),$
and the parameterized uncertainty $\tilde{\psi}_i(\chi_i) = \phi_i(\chi_i)^\top \theta$, with $\theta = (-1,1,0.5)$ and $\phi_i(\chi_i) =\left(\sin\left(\chi_{1,i}\right),~ |\chi_{2,i}|\chi_{2,i},~ e^{\chi_{2,i}\chi_{1,i}}\right)$, for all $i\in\mathcal{V}$. For the MRAC controllers, we consider a second order reference model with natural frequency and damping ratio equal to $1$, a state-feedback gain $K=(1,1)$, and a feed-forward term $u_{f_i}\left(\chi_r\right) = -\frac{2i}{2i-1}\left(\mathbf{1}_2^\top \chi_r - r\right)$, for all $i\in\mathcal{V}$. 
Each agent records two measurements of $\tilde{\psi}_i$ and $\chi_i$ at times $t_k\in\{0,1.5\}$. The corresponding data matrices $\Delta_i$ are not individually rich, which precludes the direct application of standard CL techniques \cite{chowdhary2010concurrent}. However, the collective data satisfies the CSR condition \eqref{assumption:CSR} with $\alpha = 0.9$. To regulate the state $\chi_i$ to zero, we choose $r=0,~k_{r} =1,~k_{t} = 0,~k_{c}=0.1$, $k_{a} = 3$, $T_0=0.1$, and $T=5$ to avoid instability in the estimation due to the asymmetry of the graph. We let each agent implement an MRAC controller interconnected with the hybrid dynamics $\mathcal{H}$ and show the resulting trajectories in Figure \ref{fig:hybridMRAC}. As observed, the DMCL algorithm with restart yields better transient performance compared to traditional first-order cooperative approaches without momentum \cite{poveda2021data}. Note that these results are obtained using decentralized  \emph{recorded} (i.e., batch) data, as opposed to real-time PE data. The latter might require potentially extreme transient excursions of some states whenever the parameter estimation is accelerated, which is a well-known challenge in real-time adaptive control, see \cite{zang1990transient}.
%
%
\begin{figure*}[t]
\centering
\includegraphics[width=0.97\linewidth]{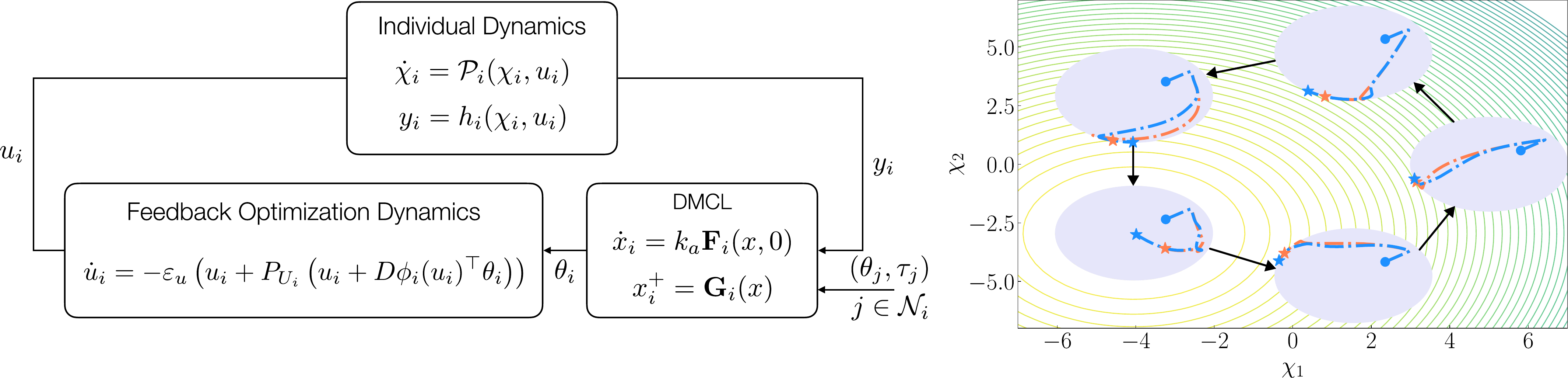}
\vspace*{-0.4cm}
\caption{Left: Scheme of the $i^{th}$ agent's dynamics in the data-enabled hybrid cooperative feedback optimization dynamics. Right: Trajectories of the vehicles. The arrows represent the edges of $\mathcal{G}$. The final positions of the vehicles are represented by stars.}
\label{fig:optimizationDynA}
\end{figure*}
\subsection{Data-Enabled Hybrid Cooperative Feedback Optimization}
Consider a multi-agent system with dynamics
\begin{equation}\label{plantdynamicsi}
\dot{\chi}_i=\mathcal{P}_i(\chi_i,u_i),~~~y_i=h_i(\chi_i,u_i),
\end{equation}
where $\chi_i\in\mathbb{R}^{n}$ is the state, $u_i\in U_i\subset \mathbb{R}$ is the input,
and $y_i\in\mathbb{R}$ is the output. The set $U_i$ is assumed to be compact and convex for all $i \in \mathcal{V}$. We consider the setting where agents seek to cooperatively find, in real-time and in a model-free manner, an optimal input $u^*$ that maximizes their individual outputs at \emph{steady state}. This scenario describes a classic data-enabled model-free feedback optimization or extremum-seeking problem  \cite{poveda2021data}. To guarantee that this problem is well-posed, the function $\mathcal{P}\coloneqq \mathcal{P}_1\times \mathcal{P}_2\times \ldots \times \mathcal{P}_N$ is assumed to be globally Lipschitz in both arguments, and we also assume there exists a smooth function $u\mapsto m(u)=m_1(u_1)\times m_2(u_2)\times \ldots\times m_{N}(u_N)$, such that for each fixed $u\in U\coloneqq U_1\times U_2\times\ldots\times U_N\subset\mathbb{R}^{N}$, the system $\dot{\chi}=\mathcal{P}(\chi,u)$ renders the equilibrium point $\chi^\star=m(u)$ UGES, uniformly on $u$.  Since the function $m(\cdot)$ describes the steady-state input-to-state mapping  of \eqref{plantdynamicsi}, the optimization problem that each agent $i$ seeks to solve can be written as
\begin{equation}\label{steadystateoptimization}
~\max_{u_i\in U_i}~J_i(u_i)\coloneqq h_i(m_i(u_i),u_i),
\end{equation}
where the \emph{response maps} $J_i$ are assumed to be unknown, continuously differentiable, strongly convex, common across the network; and parametrizable as $J_i(u_i)=\phi_i(u_i)^\top \theta^\star$, for all $u_i\in U_i$, where $\phi_i$ is a known  continuous and bounded regressor. Functions that satisfy these conditions are common in source seeking problems, where a group of mobile robots seeks to cooperatively find the maximizer of a \emph{common} potential field using intensity measurements, see \cite{poveda2021data}. In the more general case, we note that, by the universal approximation property of smooth functions, the above assumption on $J$ always holds on compact sets, modulo a small residual error that is also bounded on compact sets. In this case (i.e., non-zero approximation error), our result still holds but now in a ``semi-global practical'' sense, provided that the bound on the residual approximation error is sufficiently small, a property that can always be achieved by increasing the complexity (i.e., number of basis functions, etc) of the approximator.
%
%
%
%
%
\subsubsection{Three-Time Scale Hybrid Dynamics}
To solve the model-free feedback optimization problem using recorded data that is distributed among the agents, we use a three-time scale approach. Let $u^\star=(u_1^\star,u_2^\star,\ldots,u_N^\star)$ be the vector whose entries are the solutions to the $N$ optimization problems defined in \eqref{steadystateoptimization}. To steer $u$ towards $u^*$, we consider the following feedback optimization dynamics:
\begin{equation}\label{optimizationalgorithmi}
\dot{u}_i=-\varepsilon_u u_i+\varepsilon_u P_{U_i}\left(u_i+D\phi_i(u_i)^\top \theta_i\right),~~\forall~i\in\mathcal{V},
\end{equation}
where $D\phi_i(u_i)$ is the Jacobian matrix of $\phi_i(u_i)$, 
%
the function $P_{U_i}(\cdot)$ is the Euclidean projection on the set $U_i$, $\varepsilon_u>0$ is a tunable parameter, and $\theta_i$ is the individual estimate of $\theta^\star$, which will be recursively updated using the DMCL algorithm with restart, modeled by the hybrid system $\mathcal{H}$; refer to Figure \ref{fig:optimizationDynA} for an illustration of the overall control scheme.

To study the stability properties of the closed-loop system, we modeled the overall dynamics as a three-time scale system, where the plant dynamics \eqref{plantdynamicsi} operate at a faster time scale, the DMCL dynamics with restart operate in a medium time scale, and the optimization dynamics \eqref{optimizationalgorithmi} operate at the slowest time scale. Such time scale separation can be induced by an appropriate tuning of the gains $\varepsilon_u$ in \eqref{optimizationalgorithmi} and $k_a$ in \eqref{flomapnormalized}. By the stability assumptions on the plant dynamics \eqref{plantdynamicsi}, and by using a standard converse Lyapunov theorem \cite[Thm. 4.14]{Khalil:1173048}, there exists a Lyapunov function $V_1:\mathbb{R}^{nN}\to\mathbb{R}$, and constants $c_i>0$, for $i\in\{1,2,3,4\}$, such that $c_1|\chi-m(u)|^2\leq V_1(\chi)\leq c_2|\chi-m(u)|^2$, $\langle \nabla V_1(\chi),P(\chi,u)\rangle \leq -c_3 V_1(\chi)$,~and $|\nabla V_1(\chi)|\leq c_4|\chi-m(u)|$ for all $\chi\in\mathbb{R}^n$ and $u\in U$. Similarly, by the proof of Theorem \ref{theorem:decentralized}, and since the HDS $\mathcal{H}$ satisfies the hybrid basic conditions \cite[Ch.6]{bookHDS}, there exists a \emph{quadratic} Lyapunov function $V$ that decreases exponentially fast during flows and jumps of $\mathcal{H}$, provided that the data matrices $\{\Delta_i\}_{i\in\mathcal{V}}$ are CSR. 
Additionally, since the static-map \eqref{steadystateoptimization} is convex, the optimization dynamics \eqref{optimizationalgorithmi} with $\theta_i=\theta^\star$ reduced to a projected gradient flow that renders UGES the point $u_i^\star$ via the quadratic Lyapunov function $V_2=\frac{1}{2}|u_i-u_i^*|^2$, which satisfies $\dot{V}_2\leq -\gamma_2 V_2$ \cite[Thm. 4]{gao2003exponential}. Using these individual quadratic-type Lyapunov functions, and the global Lipschitz properties of the vector fields \eqref{plantdynamicsi}, \eqref{flomapnormalized}, and \eqref{optimizationalgorithmi}, we can now use the Lyapunov function $\hat{V}=V+V_1+V_2$ to establish exponential stability of the closed-loop system by following, recursively, the exact same steps of \cite[Ch. 11.5]{Khalil:1173048}, and using sufficiently small gains $\varepsilon_u$ and $k_a$.

\subsubsection{Numerical Example}
\begin{figure}[t!]
\includegraphics[width=0.99\linewidth]{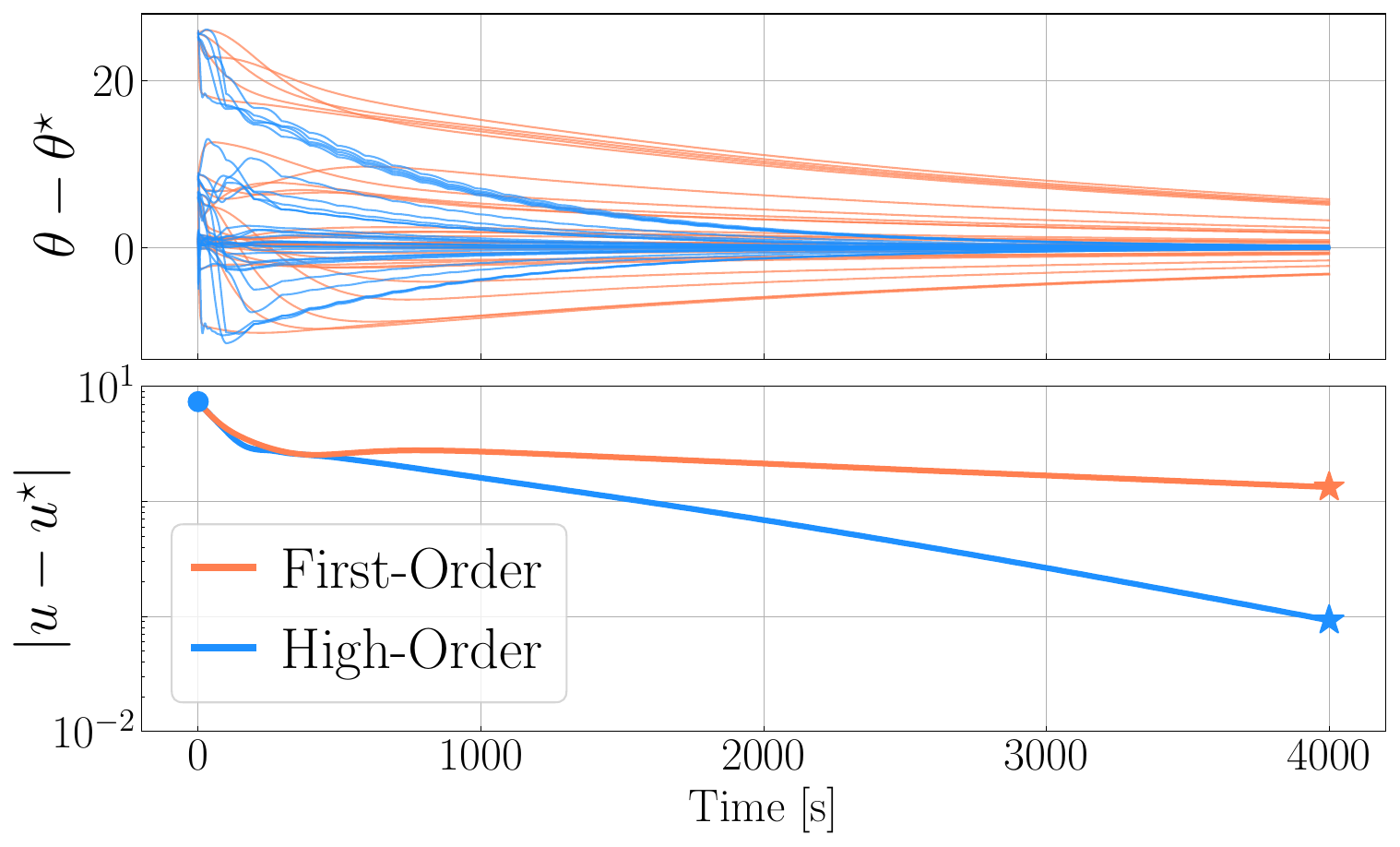}
\caption{Evolution in time of parameter (top) and control (bottom) errors.}
\label{fig:optimizationDynB}
\end{figure}

Consider a multi-vehicle system with $N=5$ vehicles, seeking to collaboratively locate the source of a potential field that is only accessible via intensity measurements. The vehicles share information via a communication graph $\mathcal{G}$ characterized again by a cycle. We assume the plant dynamics \eqref{plantdynamicsi} have the form $\mathcal{P}_i= A_i\chi_i + B_i u_i$ with matrices $A_i = -iI_{2},~B_i=iI_{2}$, and quadratic output $y_i=\chi_i^\top Q_i \chi_i + w_i^\top \chi_i + d_i$ where $Q_i= -I_{2},~w_i=(-8.1, -5.88)$, and $d_i=-25$ for all $i\in\mathcal{V}$. The sets $U_i$ are given by $U_i = \xi_i + 2\mathbb{B}$ where $\xi_i = R(2\pi i/N)(1,0)$, with $R(\alpha)$ being the standard $2\times 2$ matrix that rotates a vector by an angle $\alpha$. In this case, the steady-state input-to-output map \eqref{steadystateoptimization} reduces to $J_i(u_i) = -|u_i|^2 + w_i^\top u_i + d_i$, and each agent uses the vector of basis functions $\phi_i(u_i) = \left(u_{i,1}^2,u_{i,1},u_{i,2}^2,u_{i,2},u_{i,1}u_{i,2},1\right)$, where the parameter $\theta^\star =  (-1, -8.09, -1, -5.88, 0, -25)$ is assumed to be unknown. To implement the DMCL dynamics with restart, each agent has access to only two points of data $\{u_{i,k},y_{i,k}\}_{k=1}^2$. 
%
%
In this way, while the collected data is not persistently exciting for each agent, the collective data satisfies Assumption \ref{assumption:CSR} with $\alpha = 0.75$. Using these data and the parameters $k_{r} =1,~k_{t} = 0,~k_{c}=0.1$, $k_{a} = 0.1$, $\varepsilon_u=0.01$, $T_0=0.1$, and $T=5$, we simulate the closed-loop system comprised of the plant dynamics, the optimization dynamics in \eqref{optimizationalgorithmi}, and the hybrid dynamics $\mathcal{H}$. Figure \ref{fig:optimizationDynA} shows the resulting trajectories of the vehicles, converging to the maximizers of $J_i$ in $U_i$. Figure \ref{fig:optimizationDynB} shows the evolution in time of the parameter estimation error and the control signals. It can be observed that, given the low richness of the collected data (small $\alpha$), the proposed decentralized concurrent learning algorithm with momentum achieves faster convergence compared to the first-order cooperative estimation approach of \cite{poveda2021data}.  
\section{Proofs}\label{sec:proofs}
In this section, we present the proofs and analyses of our main results.
\subsection{Proof of Proposition \ref{lemma:definitionQ}}
For the purpose of clarity, we divide the proof of Proposition \ref{lemma:definitionQ} into multiple lemmas. 

\vspace{0.1cm}
\begin{lemma}
Suppose that Assumption \ref{assumption:CSR} holds; then, there exists a unit vector $q\in\mathbb{R}^N$ such that items (a), (b), and (c) of Proposition \ref{lemma:definitionQ} hold. \QEDB 
\end{lemma}
\vspace{0.1cm}

\textbf{Proof:} Items (a)-(b) follow directly by \cite[Prop. 1]{zhang2015constructing}. To show item (c), we use the expressions in \eqref{decomposition1eq} and \eqref{eq:sigma_omega}, and by direct substitution we obtain:
\begin{align*}
\mathbf{\Sigma}+\mathbf{\Omega}&=k_{r}\mathbf{Q}\mathbf{\Delta}\!+\! \frac{k_\text{c}}{2}\left(\textbf{Q}\textbf{L}\!+\!\textbf{L}^\top\textbf{Q}\right)+\frac{k_\text{c}}{2}\left(\textbf{Q}\textbf{L}\!-\!\textbf{L}^\top\textbf{Q}\right)\\
&=k_{r}\mathbf{Q}\mathbf{\Delta}+k_\text{c}\textbf{Q}\textbf{L}.
\end{align*}
Applying a left-multiplication by $\mathbf{Q}^{-1}$ and a right-multiplication by $\tilde{\theta}$ leads to
\begin{equation}
\mathbf{Q}^{-1}\left(\mathbf{\Sigma}+\mathbf{\Omega}\right)\tilde{\theta}=k_{r}\mathbf{\Delta}\tilde{\theta}+k_\text{c}\textbf{L}\tilde{\theta},
\end{equation}
and since $\tilde{\theta}=\theta-\mathbf{1}_N\otimes\theta^{\star}$, and $\textbf{L}(\mathbf{1}_N\otimes\theta^{\star})=0$, we obtain:
\begin{align*}
\mathbf{Q}^{-1}\left(\mathbf{\Sigma}+\mathbf{\Omega}\right)\tilde{\theta}=k_{r}\mathbf{\Delta}\tilde{\theta}+k_\text{c}\textbf{L}\theta.
\end{align*}
Finally, we show that $\mathbf{\Phi}(\theta, 0)=\mathbf{\Delta}\tilde{\theta}$. Indeed, since $\mathbf{\Phi}(\theta, 0)=(\Phi_1(\theta_1,0),\ldots,\Phi_N(\theta_N,0))$ and $\Phi_i(\theta_1,0)$ is given by \eqref{eq:learningMaps2}, we have:
\begin{align*}
\Phi_i(\theta_i,0)&=\sum_{k=1}^{\bar{k}_i}\phi_i(t_{i,k})\left(\phi_i(t_{i,k})^\top \theta_i- \phi_i(t_{i,k})^\top \theta^{\star}\right)\\
=\sum_{k=1}^{\bar{k}_i}&\phi_i(t_{i,k})\phi_i(t_{i,k})^\top\left(\theta_i-\theta^{\star}\right)=\Delta_i\tilde{\theta}_i,~~\forall~i\in\mathcal{V},
\end{align*}
which implies $\mathbf{\Phi}(\theta, 0)=\textbf{diag}(\{\Delta_1,\ldots,\Delta_N\})\tilde{\theta}=\mathbf{\Delta}\tilde{\theta}$. \hfill $\blacksquare$

\vspace{0.1cm}
\begin{lemma}
There exists $\chi\in\mathcal{K}_{\infty}$ such that \eqref{eq:omegabound} holds. \QEDB 
\end{lemma}

\vspace{0.1cm}
\textbf{Proof:} Consider the following matrix:
\begin{equation}
\mathbf{W}(t)\coloneqq  \Bigl[ \mathbf{\Omega} + k_{t} \mathbf{Q}  \mathbf{A}(t) \biggr]   \biggl[ \mathbf{\Omega} + k_{t} \mathbf{Q}  \mathbf{A}(t) \Bigr]^\top,
\end{equation}
and recall that for any symmetric matrix $A\in\mathbb{R}^{n\times n}$ we have $A \preceq \lambda_{\max}(A)I_n$  \cite[Cor. 10.4.2]{bernstein2009matrix} and $\lambda_{\max}(A)\le \sigma_{\max}(A)=\lVert A\rVert$ \cite[Fact 7.12.9]{bernstein2009matrix}, where $\lambda_{\max}(A)$ and $\sigma_{\max}(A)$ are the maximum eigenvalue and the maximum singular value of $A$, respectively. By using these facts, together with the sub-multiplicativity of the  matrix norm, we obtain that:
\begin{align}
    \mathbf{W}(t)&= 
    \mathbf{\Omega} \mathbf{\Omega}^\top + k_t\left(\mathbf{\Omega}\mathbf{A}(t)^\top\mathbf{Q} +  \mathbf{Q}\mathbf{A}(t)\mathbf{\Omega}^\top\right)\notag\\
    &~~+k_t^2\mathbf{Q}\mathbf{A}(t)\mathbf{Q}\mathbf{A}(t)^\top\notag\\
    \preceq 
    &~\big(\overline{\sigma}_{\mathbf{\Omega}}^2 + 2\overline{\sigma}_{\mathbf{\Omega}}\overline{\sigma}_{\mathbf{Q}}\lVert \mathbf{A}(t) \rVert k_{t}+ \overline{\sigma}_{\mathbf{Q}}^2\lVert \mathbf{A}(t) \rVert^2 k_{t}^2\big) I_{Nn}.\label{lemma:chiOmega:ineq}
\end{align}
Since $\phi_i(\cdot)$ is uniformly bounded,  there exists $\overline{\phi}>0$ such that $\phi_i(t)<\overline{\phi}$ for all $i\in\mathcal{V}$ and all $t\in\mathbb{R}$. Combining this fact with the diagonal structure of $\mathbf{A}(t)$ leads to $\lVert\mathbf{A}(t)\rVert \le (\overline{\phi})^2$. By using this bound in \eqref{lemma:chiOmega:ineq} we obtain:
\begin{align*}
    \mathbf{W}(t) \preceq \left(\overline{\sigma}_{\mathbf{\Omega}}^2 + 2\overline{\sigma}_{\mathbf{\Omega}}\overline{\sigma}_{\mathbf{Q}}\overline{\phi}^2 k_{t}+ \overline{\sigma}_{\mathbf{Q}}^2\overline{\phi}^4 k_{t}^2\right)I_{Nn}.
\end{align*}
The result follows using $\chi(k_t)\coloneqq \sqrt{2\overline{\sigma}_{\mathbf{\Omega}}\overline{\sigma}_{\mathbf{Q}}\overline{\phi}^2 k_{t}+ \overline{\sigma}_{\mathbf{Q}}^2\overline{\phi}^4 k_{t}^2}$, which is clearly a class-$\mathcal{K}_{\infty}$ function. \strut\hfill$\blacksquare$

\subsection{Proof of Proposition \ref{proposition2}}
We divide the proof into two lemmas:
\vspace{0.1cm}
\begin{lemma}\label{propostion:positiveDefiniteSigma}
Under Assumption \ref{assumption:CSR}, item (d) of Proposition \ref{lemma:definitionQ} holds, i.e., $\mathbf{\Sigma}$ is positive definite. \QEDB 
\end{lemma}

\vspace{0.1cm}
\textbf{Proof:} We present the proof step-by-step.

\vspace{0.1cm}
\noindent 
(a) First, note that $\mathbf{Q}\mathbf{\Delta}=\mathbf{\Delta}\mathbf{Q}$ since $\mathbf{Q} = \mathcal{Q}\otimes I_n = \textbf{diag}\left(\{q_1I_n,\dots, q_NI_n\}\right)$, $\mathbf{\Delta}= \text{diag}\left(\{\Delta_1, \dots, \Delta_N\}\right)$, with $\Delta_i \coloneqq \sum_{k=1}^{\overline{k}_i}\phi(s_{i,k})\phi(s_{i,k})^\top\in\mathbb{R}^{n\times n}$, and $q_iI_n\Delta_i= \Delta_iq_iI_n$ trivially. Then, since $\mathbf{Q}\succ 0$ and $\mathbf{\Delta}\succeq 0$ it follows that $\mathbf{Q}\mathbf{\Delta}\succeq 0$.

\vspace{0.1cm}
\noindent 
(b) Let the eigenvalues of the matrix $\mathcal{L}^\top \mathcal{Q} + \mathcal{Q}\mathcal{L}$ be organized as $0{=}\lambda_1{<}\lambda_2{\le}\cdots {\le} \lambda_N,$
and let $v_i\in\mathbb{R}^N$ be the eigenvector that corresponds to  the eigenvalue $\lambda_i$ and satisfies $|v_i|=1$. It follows that $v_1=\frac{1}{\sqrt{N}}\mathbf{1}_N$. 

\vspace{0.1cm}
\noindent 
(c) Let $\mathbf{M}:=\mathbf{L}^\top \mathbf{Q}+\mathbf{Q}\mathbf{L}$, and let
\begin{align*}
\mathbf{E}&~{\coloneqq}\frac{1}{\sqrt{N}}\left[ \mathbf{1}_N\otimes e_1, \cdots, \mathbf{1}_N \otimes e_n\right]\\
\mathbf{U}&\coloneqq\left[v_2 \otimes e_1 ,\cdots,  v_2\otimes  e_n,\cdots,v_N \otimes e_1,\cdots, v_N \otimes e_n\right]
\end{align*}
 %
%
where the vectors $e_i$ denote the standard basis in $\mathbb{R}^n$. Note that the matrices $\mathbf{E}\in\mathbb{R}^{Nn\times n}$ and $\mathbf{U}\in \mathbb{R}^{Nn\times (N-1)n}$ characterize the null space and the range space of $\mathbf{M}$, respectively. 

\vspace{0.1cm}
\noindent 
(d) Let $\hat{\mathbf{x}} \in \mathbb{R}^{Nn}$ be a unit vector, which we can write as 
\begin{equation}\label{xhatconstruction}
\hat{\mathbf{x}}=\mathbf{E}\mathbf{b}+\mathbf{U}\mathbf{c}
\end{equation}
where
$\mathbf{b}  \in \mathbb{R}^{n}$ and $\mathbf{c} \in \mathbb{R}^{(N-1)n}$ satisfy $|\mathbf{b}|^2+|\mathbf{c}|^2=1.$ 

\vspace{0.1cm}
\noindent 
(e) Since $\mathbf{E}$ can be written as $\mathbf{E}=\frac{1}{\sqrt{N}}\mathbf{1}_N\otimes I_n$, and $\mathbf{Q}\mathbf{\Delta}=\textbf{diag}\left(\{q_1\Delta_1,\ldots,q_N\Delta_N\}\right)$, we have that
\begin{equation*}
\mathbf{Q}\mathbf{\Delta}\mathbf{E}=\frac{1}{\sqrt{N}}\left(\begin{array}{c}
q_1\Delta_1\\
\vdots\\
q_N\Delta_N\\
\end{array}\right),
\end{equation*}
which leads to 
\begin{equation}\label{importantcondition001}
\mathbf{E}^\top \mathbf{Q}\mathbf{\Delta} \mathbf{E}=\frac{1}{N}\sum_{i=1}^Nq_i\Delta_i.
\end{equation}
Using \eqref{xhatconstruction} and \eqref{importantcondition001}, we obtain
\begin{align*}
\hat{\mathbf{x}}^\top\mathbf{Q}   \mathbf{\Delta} \hat{\mathbf{x}}&\geq \mathbf{b}^\top\mathbf{E}^\top \Delta \mathbf{E} \mathbf{b}+2\mathbf{b}^\top \mathbf{E}^\top\mathbf{Q}\mathbf{\Delta}\mathbf{U}\mathbf{c}\\
&\geq \underline{\sigma}_{\mathbf{Q} }\bar{\alpha} |\mathbf{b}|^2+
		2\mathbf{b}^\top \mathbf{E}^\top \mathbf{Q}  \mathbf{\Delta} \mathbf{U}\mathbf{c},
\end{align*}
where $\bar{\alpha}:=\alpha/N$, 
%
%
$\alpha$ is given by Assumption \ref{assumption:CSR}, $\underline{\sigma}_{\mathbf{Q} }$ and $\overline{\sigma}_{\mathbf{Q} }$ are defined in \eqref{eq:Q_matrix}, and $\overline{\sigma}_{\mathbf{\Delta}} \coloneqq \lvert \mathbf{\Delta}\rvert$. Moreover, since $|2\mathbf{b}^\top \mathbf{E}^\top \mathbf{Q}  \mathbf{\Delta} \mathbf{U}\mathbf{c}|\leq 2|\mathbf{b}||\mathbf{c}|\overline{\sigma}_{\mathbf{Q} }  \overline{\sigma}_{\mathbf{\Delta}}$, and using $|\mathbf{c}|=\sqrt{1-|\mathbf{b}|^2}$, we obtain:
\begin{align}\label{xi1one}
\!\!\hat{\mathbf{x}}^\top\mathbf{Q}   \mathbf{\Delta} \hat{\mathbf{x}}&\geq \underline{\sigma}_{\mathbf{Q} }\bar{\alpha} |\mathbf{b}|^2-
		2\overline{\sigma}_{\mathbf{Q} }  \overline{\sigma}_{\mathbf{\Delta}}  |\mathbf{b}|\sqrt{1-|\mathbf{b}|^2}=:\xi_1(\mathbf{b}).
\end{align}

\noindent 
(f) On the other hand, we have that 
\begin{equation}\label{xi_2two}
\hat{\mathbf{x}}^\top  \mathbf{M} \hat{\mathbf{x}}\geq \lambda_2 |\mathbf{c}|^2=\lambda_2(1-|\mathbf{b}|^2)=:\xi_2(\mathbf{b}).
\end{equation}
%
%
%
Since by the construction of $\mathbf{\Sigma}$ in \eqref{eq:sigma_omega} we have $\hat{\mathbf{x}}^\top  \mathbf{\Sigma} \hat{\mathbf{x}}= k_{r}\hat{\mathbf{x}}^\top \mathbf{Q}  \mathbf{\Delta} \hat{\mathbf{x}}+\frac{k_c}{2} \hat{\mathbf{x}}^\top  \mathbf{M} \hat{\mathbf{x}}$, the above bounds imply that  $\mathbf{\Sigma}\succeq \underline{\sigma}_{\mathbf{\Sigma}} I_{Nn}$, where 
\begin{equation}\label{lemma:crucial:lowerbound}
	\underline{\sigma}_{\mathbf{\Sigma}}~~\geq ~\min_{0\leq \nu\leq 1} \max \left \{k_r\xi_1(\nu),\frac{k_c}{2}\xi_2(\nu)\!\right\},
\end{equation} 
with $\xi_1$ given by \eqref{xi1one} and $\xi_2$ given by \eqref{xi_2two}.

\vspace{0.1cm}
\noindent 
(g) Next, we study \eqref{lemma:crucial:lowerbound} and we show that this lower bound is indeed positive. Since, by item (a), $\mathbf{Q}\mathbf{\Delta}\succeq 0$, without loss of generality we can assume that the first term in the brackets in \eqref{lemma:crucial:lowerbound} is non-negative. Indeed, suppose by contradiction that such term is negative. Then, since $\mathbf{Q}\mathbf{\Delta}\succeq 0$, we can take $\xi(\mathbf{b})$ as a non-negative lower bound for $\hat{\mathbf{x}}^\top\mathbf{\Sigma}\hat{\mathbf{x}}$, and since $\xi(\mathbf{b})=0$ only if $\mathbf{b}=1$, we obtain that for such $\mathbf{b}$ the first term in the brackets is indeed positive.  

\vspace{0.1cm}
\noindent 
(h) To get a closed form of the expression in \eqref{lemma:crucial:lowerbound}, let $\nu=\sin(\theta),~\theta\in [0,\pi/2]$. In the $\theta$ variable, \eqref{lemma:crucial:lowerbound} becomes:
{\small
\begin{align*}
		\min_{0\leq \theta \leq \frac{\pi}{2}}  \max \Big \{k_1(1 {-} \cos(2\theta)) {-}
		k_2 \sin(2\theta),k_3(1+\cos(2\theta)) \Big\},
\end{align*}
}
where the constants $k_1,k_2,k_3>0$ are given by $k_1\coloneqq  \frac{k_{r} \underline{\sigma}_{\mathbf{Q} }\bar{\alpha}}{2},~k_2\coloneqq  k_{r}
	\overline{\sigma}_{\mathbf{\Delta}}  \overline{\sigma}_{\mathbf{Q} },~k_3\coloneqq  \frac{k_c}{4}\lambda_2$. Further simplifying, we obtain
\begin{align}
		&\min_{0\leq \theta \leq \frac{\pi}{2}}  \max \Bigg \{  k_1-\sqrt{k_1^2+k_2^2}\sin\left(2\theta+\tan^{-1}\left(\frac{k_1}{k_2}\right)\right), \notag\\&~~~~~~~~~~~~~~~~~~~~~k_3+k_3\cos(2\theta)) \Bigg\}\notag\\
		&~~~~\coloneqq  \min_{0 \le \theta \le \frac{\pi}{2}} \max \Big \{ g_1(\theta),~g_2(\theta) \Big\}.\label{minmaxproblem}
\end{align}
We argue that the intersection point  $\theta^* \in [0,\frac{\pi}{2}]$ of the trigonometric curves $g_1(\theta),~g_2(\theta)$  solves the min-max problem \eqref{minmaxproblem}.

\begin{figure}[t!]
    \centering
    \includegraphics[width=0.9\linewidth]{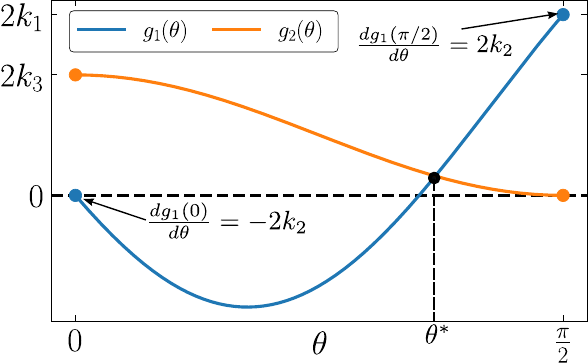}
    \caption{Illustration of step (i) in the proof of Proposition \ref{propostion:positiveDefiniteSigma}. }
    \label{fig:lemma_fig_PD}
\end{figure}

\vspace{0.1cm}
\noindent 
(i) To establish the existence of such $\theta^* \in [0,\frac{\pi}{2}]$, we use the following facts: (i) $k_1, k_2,k_3>0$, (ii) $g_1(0)=0$, $g_1(\frac{\pi}{2})=2k_1>0$, $\frac{d g_1(\theta)}{d \theta}={-2\sqrt{k_1^2+k_2^2}\cos(2\theta+\tan^{-1}(\frac{k_1}{k_2}))}$, (ii) $g_2(0)=2k_3>0$, $g_2(\frac{\pi}{2})=0$, $\frac{d g_2(\theta)}{d \theta}=-2k_3\sin(2\theta)$.  Since $g_1$ and $g_2$ are continuous functions, the previous conditions imply the existence of a point $\theta^*$ such that $g_1(\theta^*)=g_2(\theta^*)$. Moreover, since $g_2$ is decreasing on $[0, \frac{\pi}{2}]$ with $g_2(0)>0$, $g_2(\frac{\pi}{2})=0$, $g_1(0)=0$, $g_1(\frac{\pi}{2})>0$, and $\frac{d g_1(\theta)}{d \theta}={-2\sqrt{k_1^2+k_2^2}\cos(2\theta+\tan^{-1}(\frac{k_1}{k_2}))}$, it follows that the intersection point $\theta^*$ is in fact the minimum of the point-wise maximum of $g_1(\theta)$ and $g_2(\theta)$. See Figure \ref{fig:lemma_fig_PD} for an illustration of this step.

\vspace{0.1cm}
\noindent 
(j) By computing the intersection point $\theta^*$, we obtain
{\small
\begin{align*}
	\theta^*=\frac{1}{2}\left[\cos^{\!-1}\!\!\left(\frac{k_1-k_3}{\sqrt{(k_1+k_3)^2+k_2^2}}\right)+\tan^{\!-1}\!\!\left(\frac{k_2}{k_1+k_3}\right)\right].
\end{align*}
}
Substituting the values of $k_1$, $k_2$ and $k_3$, establishes the existence of a positive lower bound on the constant $\underline{\sigma}_{\mathbf{\Sigma}}$ that satisfies $\mathbf{\Sigma} \succeq \underline{\sigma}_{\mathbf{\Sigma}} I_{Nn}$, given by 
\begin{equation}
\underline{\sigma}_{\mathbf{\Sigma}}\ge \frac{k_c \lambda_2}{4}\bigl[1+\cos( \theta^*)\bigr],
\end{equation}
%
%
where $\theta^*=\theta_1^*+\theta_2^*$, with 
\begin{equation*}
\theta_1^*=\cos^{-1}\left(\frac{2k_{r}\underline{\sigma}_{\mathbf{Q} }\alpha-k_c\lambda_2}{\sqrt{(2k_{r}\underline{\sigma}_{\mathbf{Q} }\alpha+k_c\lambda_2)^2+16 k_{r}^2  \overline{\sigma}_{\mathbf{\Delta}} ^2 \overline{\sigma}_{\mathbf{Q}}^2}}\right)
\end{equation*}
and $\theta_2^*=\tan^{-1}\left(\frac{4  \overline{\sigma}_{\mathbf{\Delta}}  k_{r} \overline{\sigma}_{\mathbf{Q}}}{2\underline{\sigma}_{\mathbf{Q} }k_{r}\alpha +k_c\lambda_2}\right)$. Note that $\cos(\theta^*)\in[0,1]$ since $\theta^*\in [0,\pi/2]$, which implies that $~\underline{\sigma}_{\mathbf{\Sigma}}>0$. \hfill $\blacksquare$
%
%

%

\vspace{0.1cm}
\begin{lemma}\label{auxlemma00}
Let $\lambda_N$ be the largest eigenvalue of $\mathcal{L}^\top \mathcal{Q} + \mathcal{Q}\mathcal{L}$. Then, under Assumption \ref{assumption:CSR}, the matrix $\mathbf{\Sigma}$ satisfies 
\begin{equation}\label{upperboundsigma}
\mathbf{\Sigma}  \preceq \left(k_{r} \overline{\sigma}_{\mathbf{Q} } \overline{\sigma}_{\mathbf{\Delta}} + \frac{k_c}{2} \lambda_N\right) I_{Nn}.
\end{equation}
\end{lemma}

\textbf{Proof:} By the definition of $\mathbf{\Delta}$ and $\overline{\sigma}_{\mathbf{\Delta}}$, the term $\mathbf{Q} \mathbf{\Delta}$ satisfies: 
$
\mathbf{Q} \mathbf{\Delta} \preceq  \overline{\sigma}_{\mathbf{Q} } \overline{\sigma}_{\mathbf{\Delta}}  I_{Nn}$. By the definition of $\lambda_N$ and the fact that $\mathbf{Q}\mathbf{L}= \mathcal{Q}\mathcal{L}\otimes I_n$ by the properties of the Kronecker product, it follows that $ \mathbf{L}^\top \mathbf{Q} +\mathbf{Q} \mathbf{L}  \preceq  \lambda_N I_{Nn}$.  Note that $\lambda_N>0$ since, as stated in the proof of Lemma \ref{lemma:definitionQ}, $\mathcal{Q}\mathcal{L} + \mathcal{L}^\top \mathcal{Q}$ is a nonzero and symmetric $M$-matrix.  Combining these arguments we obtain \eqref{upperboundsigma}.

\subsection{Proof of Theorem \ref{theorem:centralized}}
We follow a (hybrid) Lyapunov-based approach to study the HDS $\mathcal{H}_c$ with input $u$, in the error coordinates $$\tilde{y}_c=(\tilde{x}_c,s)\coloneqq ((\tilde{\theta},\tilde{p},\tau_c),s),$$ where $\tilde{\theta} = \theta - \mathbf{1}_N\otimes\theta^\star$, $\tilde{x}_c=(\tilde{\theta},\tilde{p},\tau_c)$, and $\tilde{p}=p-\mathbf{1}_N\otimes \theta^\star$. In these new coordinates, the HDS with input $u$ becomes $$\tilde{\mathcal{H}}_c = (\mathbf{C}_c\times\mathbb{R}_{\ge0}, \tilde{\mathbf{F}}_{c}, \mathbf{D}_c\times\mathbb{R}_{\ge0}, \mathbf{G}_c),$$ where $\tilde{\mathbf{F}}_{c}(\tilde{y}_c,u)\coloneqq \hat{\mathbf{F}}_c(\tilde{y}_c,u)\times [0,\omega]\times \{1\}$, with $ \hat{\mathbf{F}}_c$ given by \eqref{errormapping}. For this system, we will study stability properties with respec to the set $\tilde{\mathcal{A}}_c\times \mathbb{R}_{\ge 0}$, where 
\begin{equation}\label{setorigin}
\tilde{\mathcal{A}}_c\coloneqq \{0\}\times \{0\}\times [T_0,T].
\end{equation}
\subsubsection{Proof of Theorem \ref{theorem:centralized}-(a)}
We establish item (a) of Theorem \ref{theorem:centralized} via a sequence of lemmas. The following lemma follows directly from the uniform boundedness assumption on the regressors $\phi$ and the definition of $\mathbf{U}$ in \eqref{defU}.

\vspace{0.1cm}
\begin{lemma}\label{boundedregressorslemma}
There exist $\overline{\phi}>0$ such that $|\mathbf{U}(s)|\le \overline{\phi}|u|$ for all $s\geq0$. \QEDB
\end{lemma}

\vspace{0.1cm}

Next, we consider the Lyapunov function
\begin{equation}\label{eq:lyap_fcn}
		V(\tilde{y}_c)\coloneqq  \frac{|\tilde{p}-\tilde{\theta}|^2_{\mathbf{Q}}}{4}+\frac{|\tilde{p}|^2_{\mathbf{Q} }}{4}+ \tau_c^2 \frac{|\tilde{\theta}|^2_\mathbf{\Sigma}}{2}.
\end{equation} 
and we study its behavior during flows and jumps of $\tilde{\mathcal{H}}_c$. 
%
%
%
and present a lemma and two auxiliary propositions.
%
\begin{lemma}\label{proof:centralized:lemma:lyapunov_function_bounds}
There exist constants $\overline{c}>\underline{c}>0$ such that $$\underline{c}|\tilde{y}_c|^2_{\tilde{\mathcal{A}}_c\times\mathbb{R}_{\geq0}} \leq V(\tilde{y}_c) \leq \overline{c}|\tilde{y}_c|^2_{\tilde{\mathcal{A}}_c\times\mathbb{R}_{\geq0}},$$ for all $\tilde{y}_c\in (\mathbf{C}_c\cup \mathbf{D}_c)\times \mathbb{R}_{\ge0}$. \QEDB
\end{lemma}

\vspace{0.1cm}
\textbf{Proof:} Since, by the definition of $\tilde{\mathcal{H}}_c$, we always have $s\in\mathbb{R}_{\geq0}$, we just need to study $|\tilde{x}_c|_{\tilde{\mathcal{A}}_c}$. To establish the lower bound, and using the definition of the norm $|\cdot|_{\mathbf{P}}$, and since $\tau_c\geq T_{0}$ for all $\tilde{x}_c\in \mathbf{C}_c\cup \mathbf{D}_c$, we directly obtain that $|\tilde{p}|^2_{\mathbf{Q} } \ge \underline{\sigma}_{\mathbf{Q} } |\tilde{p}|^2$ and $\tau^2 |\tilde{\theta}|^2_\mathbf{\Sigma} \ge  \underline{\sigma}_{\mathbf{\Sigma}}T_{0}^2 |\tilde{\theta}|^2$. Therefore, $V(\tilde{y}_c)\ge\underline{c}|\tilde{x}_c|_{\tilde{\mathcal{A}}_c}^2$, where $\underline{c}:=\frac{1}{4}\min\left\{\underline{\sigma}_{\mathbf{Q} }, \, 2 \underline{\sigma}_{\mathbf{\Sigma}} T_{0}^2\right\}$.
%
%
%
To establish the upper bound, we use \eqref{eq:Q_matrix} together with the fact that $\tau\le T$ to obtain that $V(\tilde{x}_c,s)\le  \frac{1}{4} (2\overline{\sigma}_{\mathbf{Q} }|\tilde{\theta}|^2 + 3\overline{\sigma}_{\mathbf{Q} }|\tilde{p}|^2 + 2T^2|\tilde{\theta}|^2_\mathbf{\Sigma})$, where we also used the fact that $|\tilde{p}-\tilde{\theta}|^2\le 2(|\tilde{\theta}|^2+|\tilde{p}|^2)$.
%
Using Lemma \ref{auxlemma00}, we obtain $V(\tilde{y}_c)\leq \overline{c}|\tilde{x}_c|_{\tilde{\mathcal{A}}}^2$, with $\overline{c}:=\frac{1}{4}\max \left\{3\overline{\sigma}_{\mathbf{Q} }, \,  T^2(2k_{r} \overline{\sigma}_{\mathbf{Q} } \overline{\sigma}_{\mathbf{\Delta}}  + k_c\lambda_N) + 2\overline{\sigma}_{\mathbf{Q} }\right\}$. \hfill $\blacksquare$
%
%

\vspace{0.1cm}
\begin{lemma}\label{proof:centralized:proposition:flows}
Suppose that $T<\overline{\mathbf{T}}$; then, there exists $\varrho>0$ and $\gamma>0$ such that $V$ satisfies $\dot{V}(\tilde{y}_c) \leq\! -\varrho V(\tilde{y}_c) + \gamma|u|^2,$ for all $\tilde{y}_c\in\mathbf{C}_c\times \mathbb{R}_{\ge 0}$. \hfill\QEDB
%
%
\end{lemma}

\vspace{0.1cm}
\textbf{Proof:} By direct computation, we have:
\begin{align}
\dot{V}(\tilde{y}_c)&=- \tau_c \begin{pmatrix}(\tilde{p}-\tilde{\theta})^\top&\tilde{\theta}^\top
		\end{pmatrix}\mathbf{V}_w(\tau_c,s)\begin{pmatrix}\tilde{p}-\tilde{\theta}\\\tilde{\theta}
		\end{pmatrix}\notag\\
      &\quad +\tau_c(2\tilde{p}-\tilde{\theta})^\top\mathbf{Q}\mathbf{U}(s),\label{eq:lyap_bound2}
\end{align}
where 
\begin{equation*}
\mathbf{V}_w(\tau_c,s)\coloneqq \begin{pmatrix}\frac{\mathbf{Q} }{ \tau_c^2} &  \hat{\mathbf{\Omega}}(s)\\  \hat{\mathbf{\Omega}}(s)^\top&  (1-w)\mathbf{\Sigma} + k_t\mathbf{Q}\mathbf{A}(s)\end{pmatrix}, 
\end{equation*}
for all $w\in[0,\omega]$, where
\begin{equation*}
\hat{\mathbf{\Omega}}(s) \coloneqq \mathbf{\Omega}+k_t\mathbf{Q}\mathbf{A}(s),
\end{equation*}
and where we used the fact that $x_1^\top \mathbf{\Omega}x_1 =0$. Using the definitions of $\underline{\sigma}_{\mathbf{Q}}$, $\underline{\sigma}_{\mathbf{\Sigma}}$, and $\overline{\sigma}_{\mathbf{\Omega}}$, and Lemma \ref{lem:Vbound} in the Appendix, it follows that $\mathbf{V}_w(\tau_c,s)\succeq \underline{\mathbf{v}}I_{Nn}$, for all $\tau_c\in[T_0,T]$, all $w\in[0,\omega]$, and all $s\in\mathbb{R}_{\ge 0}$, with 
\begin{equation}
\underline{\mathbf{v}}~{\coloneqq}\frac{(1-\omega)\underline{\sigma}_{\mathbf{\Sigma}}\underline{\sigma}_{\mathbf{Q}} - T^2   (\overline{\sigma}_{\mathbf{\Omega}}^{2}+\chi(k_t)^2)}{T^2(1-\omega)\underline{\sigma}_{\mathbf{\Sigma}} + \underline{\sigma}_{\mathbf{Q}}}>0,
\end{equation}
and $\chi\in\mathcal{K}_\infty$. Using the Cauchy-Schwartz inequality to upper-bound the last term in \eqref{eq:lyap_bound2}, and since $T_0\le \tau_c\le T$ and $|\tilde{x}_c|^2\le 3(|\tilde{p}-\tilde{\theta}|^2 + |\tilde{p}|^2)$ for all $\tilde{x}_c\in \mathbf{C}_c\cup\mathbf{D}_c$, we obtain:
\begin{align}
		\dot{V}(\tilde{y}_c)&\le -T_0 \underline{\mathbf{v}} (|\tilde{p}\!-\!\tilde{\theta}|^2+|\tilde{\theta}|^2)\notag\\
            &\qquad+ 2T(|\tilde{p}| + |\tilde{\theta}|)\|\mathbf{Q}\||\mathbf{U}(s)|\notag\\
                               &\leq -\frac{\underline{\mathbf{v}}}{3} T_0 |\tilde{x}_c|^2 + 2\sqrt{2}\overline{\sigma}_{\mathbf{Q}}\overline{\phi}T|\tilde{x}_c||u|\notag\\
                              &\leq\! -\bigg(\frac{\underline{\mathbf{v}}}{3} T_0 \!-\! \frac{1}{\epsilon}\bigg)|\tilde{x}_c|^2 \!+\! 2\epsilon\left(\overline{\sigma}_{\mathbf{Q}}\overline{\phi}T\right)^2
                               |u|^2,\label{proof:decreaseflows}
\end{align}
for all $\epsilon>0$ and all $w\in[0,\omega]$, where the last inequality follows from the fact that $ab \le \frac{1}{4\epsilon}a^2 + \epsilon b^2$ for all $\epsilon>0$. Setting $\epsilon:=\frac{3(1+\varepsilon)}{T_0\nu}$ for $\varepsilon>0$ and using the lower bound of Lemma~\ref{proof:centralized:lemma:lyapunov_function_bounds}, the expression in \eqref{proof:decreaseflows} yields:
\begin{align}
    \!\!\!\dot{V}(\tilde{y}_c)&\le\! -\frac{\varepsilon}{1\!+\!\varepsilon}\frac{\nu T_0}{3\overline{c}}V(\tilde{y}_c) + (1\!+\!\varepsilon)\frac{6}{\nu T_0}\left(\overline{\sigma}_{\mathbf{Q}}\overline{\phi}T\right)^2|u|^2.\label{proof:centralized:eq:decreaseflows}
\end{align}
The result follows by setting $\varrho\coloneqq \frac{\nu T_0}{3\overline{c}}\frac{\varepsilon}{1+\varepsilon}$ and letting $\gamma\in\mathcal{K}_\infty$ be defined as $\gamma(r) \coloneqq (1+\varepsilon)\frac{6}{\nu T_0}\left(\overline{\sigma}_{\mathbf{Q}}\overline{\phi}T\right)^2r$. \hfill $\blacksquare$
\vspace{0.1cm}
\begin{lemma}\label{proof:centralized:proposition:jumps}
Suppose that $T>\underline{\mathbf{T}}$; then, 
$$V(\tilde{y}_c^+)\le  \left(\mu_T\right)^{\underline{\eta}} V(\tilde{y}_c),~~~\forall~\tilde{y}_c\in(\mathbf{C}_c\cup\mathbf{D}_c)\times\mathbb{R}_{\ge0},$$
where $\underline{\eta}:=\min_{i\in\mathcal{V}}\eta_i$, and $\mu_T\coloneqq(\underline{\textbf{T}}/T)^2$.
%
%
\QEDB
\end{lemma}

\vspace{0.1cm}
\textbf{Proof:} Using the definition of the jump map $\mathbf{G}_c$, for all $\tilde{y}_c\in\mathbf{D}_c\times\mathbb{R}_{\ge 0}$ we have:
\begin{align}
4V(\tilde{y}_c^+)&=
|\mathbf{R}_{\eta}\tilde{\theta}{+}(I_{Nn}{-}\mathbf{R}_{\eta})\tilde{p}{-}\tilde{\theta}|_{\mathbf{Q} }^2\notag  \\
& + |\mathbf{R}_{\eta}\tilde{\theta}+(I_{Nn}- \mathbf{R}_{\eta})\tilde{p}|_{\mathbf{Q}}^2+2T_0|\tilde{\theta}|^2_\mathbf{\Sigma},\label{eq:deltaVjumps}
\end{align}
where $\mathbf{R}_\eta:=\text{diag}(\eta)\otimes I_n$. By Lemma \ref{lemma:aux:ResettingPolicies} in the Appendix, the change of $V$ during jumps, given by $\Delta V(\tilde{y}_c)\coloneqq  V(\tilde{y}^+_c)-V(\tilde{y}_c)$,  satisfies:

\vspace{-0.3cm}
\begin{small}
\begin{align*}
    4\Delta V(\tilde{y}_c)&=
  |\tilde{\theta}|^2_{\mathbf{R}_\eta\mathbf{Q}} - |\tilde{p}|^2_{\mathbf{R}_\eta\mathbf{Q}} -  |\tilde{\theta} -\tilde{p}|^2_{\mathbf{R}_\eta\mathbf{Q}}\\
  &\qquad\qquad\qquad\qquad+ 2T_0^2|\tilde{\theta}_c|_{\mathbf{\Sigma}}^2-2T^2|\tilde{\theta}|_{\mathbf{\Sigma}}^2\notag\\
%
%
&\le- \left(|\tilde{p}|^2_{\mathbf{R}_\eta\mathbf{Q}} + |\tilde{\theta} -\tilde{p}|^2_{\mathbf{R}_\eta\mathbf{Q}}\right) + \frac{\overline{\sigma}_\mathbf{Q}}{\underline{\sigma}_{\mathbf{\Sigma}}}|\tilde{\theta}|^2_{\mathbf{\Sigma}}\\
&\qquad\qquad\qquad\qquad + 2T_0^2|\tilde{\theta}|_{\mathbf{\Sigma}}^2-2T^2|\tilde{\theta}|_{\mathbf{\Sigma}}^2\\
%
%
&=-\left(|\tilde{p}|^2_{\mathbf{R}_\eta\mathbf{Q}} + |\tilde{\theta} -\tilde{p}|^2_{\mathbf{R}_\eta\mathbf{Q}}\right) - (1-\mu_T)2T^2|\tilde{\theta}|_{\mathbf{\Sigma}}^2\\
&\leq-(1-\mu_T)\left( |\tilde{p}|^2_{\mathbf{R}_\eta\mathbf{Q}} +  |\tilde{\theta} -\tilde{p}|^2_{\mathbf{R}_\eta\mathbf{Q}}+2T^2|\tilde{\theta}|_{\mathbf{\Sigma}}^2\right),
\end{align*}
\end{small}

\noindent 
where we also used the fact that $\mu_T\in(0,1)$ whenever $T>\underline{\mathbf{T}}$.
It then follows that $\Delta V(\tilde{y}_c)\le 0$ for all  $\tilde{y}_c\in \mathbf{D}_c\times\mathbb{R}_{\ge0}$.  When $\underline{\eta}=1$, the previous inequality yields $\Delta V(\tilde{y}_c)\le - (1 -\mu_T)V(\tilde{y}_c)$, wich in turn implies that $V(\tilde{y}_c)\le \mu_TV(\tilde{y}_c)$. \hfill $\blacksquare$

\vspace{0.2cm}
By the construction of the dynamics of $\tau_c$, every solution to $\mathcal{H}_c$ is guaranteed to have intervals of flow with a duration of at least $(T-T_0)/\omega$ between any two consecutive jumps. Combining this fact with Lemmas \ref{proof:centralized:lemma:lyapunov_function_bounds}, \ref{proof:centralized:proposition:flows}, and \ref{proof:centralized:proposition:jumps}, it follows that $\tilde{\mathcal{H}}_c$ renders the set $\tilde{\mathcal{A}}_c$ ISS with respect to the input $u$. The ISS property of the HDS $\mathcal{H}_c$ with respect to $|\cdot|_{\mathcal{A}_c}$ follows directly by employing the change of coordinates $\tilde{y}_c\to y_c$.

\vspace{0.2cm}
\subsubsection{Proof of Theorem \ref{theorem:centralized}-(b)}
Let the initial condition satisfy $\tilde{y}_0:=((\tilde{\theta}(0,0),\tilde{p}(0,0),\tau_c(0,0)),s(0,0))\in(\mathbf{C}_c\times \mathbf{D}_c)\times\mathbb{R}_{\ge 0}$, and let $(\tilde{y}_c,u)$ be a maximal solution pair to $\tilde{\mathcal{H}}_c$ from the initial condition $\tilde{y}_0$, and satisfying  $\dot{\tau}_c(t,j)=\omega$ for all $(t,j)\in\text{dom}(\tilde{y}_c)$. By Lemma \ref{proof:centralized:proposition:flows}, we have that $\dot{V}(\tilde{y}_c)\le -\frac{\varrho}{2}V(\tilde{y}_c)$ for all $\tilde{y}_c\in \left(\mathbf{C}_c\cup\mathbf{D}_c\right)\times\mathbb{R}_{\ge 0}$ such that $V(\tilde{y}_c)\ge \frac{2\gamma}{\varrho}|u|^2$. Let 
\begin{equation}
\mathcal{R}\coloneqq \left\{\tilde{y}_c\in\mathbb{R}^{2nN + 1}\times \mathbb{R}_{\ge 0}~:~V(\tilde{y}_c)\le \frac{2\gamma}{\varrho}|u|^2_{\infty}\right\},
\end{equation}
and let $\mathbb{T} \coloneqq \sup\{\sigma\in\mathbb{R}_{\ge 0}~:~\tilde{y}_c(\tilde{t},\tilde{j})\not\in\mathcal{R},~(\tilde{t},j)\in\text{dom}(\tilde{y}_c),~0\le \tilde{t} + \tilde{j} \le \sigma\}$. Then, letting $t_j\coloneqq \min\{t\in\mathbb{R}_{\ge0}~:~(t,j)\in\text{dom}(\tilde{y}_c)\}$ for every $j\in\mathbb{Z}_{\ge0}$, and via the comparison lemma, it follows that
$
 V(\tilde{y}_c(t,j))\le e^{-\rho(t-t_j)/2}V( \tilde{y}_c(t_{j},j)),	
$
for all $(t,j)\in\text{dom}(\tilde{y}_c)$ such that $t_j+j\le t+j\le \mathbb{T}$. On the other hand, from Lemma \ref{proof:centralized:proposition:jumps}, it follows that $V(\tilde{y}_c(t_{j+1},j+1))\le  \mu_T V(\tilde{y}_c(t_{j+1},j))$, which iterating over $j$ yields:
\begin{equation}\label{eq:prevT}
V(\tilde{y}_c(t,j)) \le e^{-\varrho t/2}
\mu_T^j V(\tilde{y}_{0}),
\end{equation}
for all $(t,j)\in\text{dom}(\tilde{y}_c)$ such that $t+j\le \mathcal{T}$ and where we have used that $t_0=0$. Since $\dot{V}(\tilde{y}_c(t,j))\le 0$ if $\tilde{y}_c(t,j)\in \mathcal{R}$ and Proposition \ref{proof:centralized:proposition:jumps} holds for all $(t_j,j)\in \text{dom}(\tilde{y})$, it follows that $\tilde{y}_{c}(t,j)\in \mathcal{R}$ for all $t+j\ge \mathcal{T}$, meaning that
\begin{equation}\label{eq:posT}
V(\tilde{y}_c(t,j)) \le \frac{2\gamma}{\varrho}|u|_\infty,
\end{equation}
for all $t+j\ge \mathcal{T}$. The bounds \eqref{eq:prevT} and \eqref{eq:posT}, together with Lemma \ref{proof:centralized:lemma:lyapunov_function_bounds} and the time-invariance of $\mathcal{H}_c$, imply that
$
   |\tilde{y}_c(t,j))|_{\tilde{\mathcal{A}}_c}^2  \le\frac{\overline{c}}{\underline{c}}
\mu_T^{j}|\tilde{y}_0|_{\tilde{\mathcal{A}}_c}^2+ \frac{2\gamma}{\varrho}|u|_{(t,j)}
$
for all $(t,j)\in\text{dom}(\tilde{y}_c)$, where we also used the fact that $e^{-\varrho t/2}\le 1$ for all $t\in \mathbb{R}_{\ge 0}$.
The bound \eqref{theorem:centralized:convergencebound}, is obtained by evaluating the above bound at the hybrid times $(t_j,j)$, noting that $|\tilde{\theta}|\le |\tilde{y}_c|_{\mathcal{A}_c}$, and via the change of coordinates $\tilde{y}_c\mapsto y_c$.

\subsection{Proof of Theorem \ref{theorem:decentralized}}
\label{proofProposition2}
The proof uses the reduction principle for hybrid systems \cite[Corollary 7.24]{goebel2012hybrid}. First, note that, by construction, $\mathcal{H}$ satisfies the hybrid basic conditions \cite[Assump. 6.5]{goebel2012hybrid}. Since the flow map $\mathbf{F}$  is globally Lipschitz in $\mathbf{C}$, the HDS does not exhibit finite escape times. To study the stability properties of the system, we first intersect the flow set $\mathbf{C}$, the jump set $\mathbf{D}$, and the values of the jump map $\mathbf{G}_d$ with a compact set $K\subset \mathbb{R}^{(2n+1)N}$. Since $\tau$ already evolves in a compact set, we take $K$ only to restrict the states $(\theta,p,s)$. The new restricted system is denoted as $\mathcal{H}_K=(\mathbf{C}\cap K,\mathbf{F},\mathbf{D}\cap K, \mathbf{G}\cap K)$. Since the dynamics of the state $\tau$ are independent of $(\theta,p)$, we can directly use \cite[Prop. 1-(a)]{javed2021scalable} to conclude that, under condition (b) of Theorem \ref{theorem:decentralized}: 1) the set $K\times \mathcal{A}_{\text{sync}}$ is UGAS for the HDS $\mathcal{H}_K$, and 2) $\tau$ converges to $\mathcal{A}_{\text{sync}}$ before the hybrid time $(2t^*,2N)$. It follows that, for all solutions $(y,s)$, and all times  $(t,j)\in\text{dom}((y,s))$ such that $t\geq 2t^*$ and $j\geq 2N$, the restricted synchronized HDS behaves as having the centralized master timer $\tau_c$ of Section \ref{section:centralized}. Next, we intersect the data of the HDS $\mathcal{H}_K$ with the set $K\times \mathcal{A}_{\text{sync}}$.  For this restricted HDS, denoted $\mathcal{H}_{K,\mathcal{A}_{\text{sync}}}$, Theorem \ref{theorem:centralized} guarantees UGES of the set $\mathcal{A}$ when $u=0$. By invoking the reduction principle of \cite[Corollary 7.24]{goebel2012hybrid}, we conclude UGES of the set $\mathcal{A}$ for the HDS $\mathcal{H}_K$. Since this system has bounded solutions, and $K$ was arbitrary large, for each compact set of initial conditions $K_0$ of system $\mathcal{H}$, we can select $K$ sufficiently large such that the restriction in $\mathcal{H}_K$ does not affect the solutions from $K_0$, obtaining UGES of $\mathcal{A}$ for the original hybrid system $\mathcal{H}$. Now, since the convergence of $\tau\in \mathbb{R}^{N}$ to $\mathcal{A}_{\text{sync}}$ occurs in finite time after which the stability properties are characterized by Theorem \ref{theorem:centralized}, we obtain that $\mathcal{A}$ is UGES for $\mathcal{H}$. \hfill $\blacksquare$
\subsection{Proof of Corollary \texorpdfstring{\ref{corollary:centralized:optimalrestarting}}{1}}\label{appendixoptimalresetsH1}
First, note that $j\ge \frac{t'}{\left(T - T_0\right)/\omega}$ for any $(t',j)\in \text{dom}(y)$. Therefore, since $\mu(T)\in (0,1)$, the bound \eqref{theorem:centralized:convergencebound} implies the following slightly looser bound when $u\equiv 0$:
\begin{equation}
    |y_j|^2_{\mathcal{A}}\leq\frac{\overline{c}}{\underline{c}}\left(\mu_T^{\frac{1}{T-T_0}}\right)^{\omega t'}|y_0|^2_{\mathcal{A}}.\label{eq:rateDeterministicLooser}
\end{equation}
where $\overline{c}$ and $\underline{c}$ come from Lemma \ref{proof:centralized:lemma:lyapunov_function_bounds}. Following similar ideas to \cite{o2015adaptive,poveda2021robust}, and using the definition of $\mu(T)$, we solve the following optimization problem to maximize the rate of contraction over any window of time $t'$:
\begin{equation*}
\min_{T\in\mathbb{R}_{>0}}\phi(T)\coloneqq\mu_T^{\frac{1}{T-T_0}}.
\end{equation*}
Computing the derivative of $\phi$ with respect to $T$, and equating to zero, we obtain: $T^*=e\sqrt{\frac{\overline{\sigma}_\mathbf{Q}}{2\underline{\sigma}_{\mathbf{\Sigma}}}+T_0^2}$, which is the unique minimizer of $\phi$. By substituting $T=T^*$ in \eqref{eq:rateDeterministicLooser}, we obtain 
\begin{equation}\label{eq:exponential}
    |y_j|_\mathcal{A}^2\leq \frac{\overline{c}}{\underline{c}}e^{-\tfrac{2\omega t'}{T^*-T_0}}|y_0|^2_{\mathcal{A}}.
\end{equation}
Thus, to have $|y_j|_{\mathcal{A}}^2\le \varepsilon$ for a given $\varepsilon>0$, it suffices to have that $ t' \ge \frac{1}{2\omega}\left(T^*-T_0\right)\log\left(\frac{1}{\varepsilon}\frac{\overline{c}}{\underline{c}}|y_0|^2_{\mathcal{A}}\right).$ Moreover, note that the right hand side of \eqref{eq:exponential} is of order $\mathcal{O}\left(e^{-\sqrt{\underline{\sigma}_{\mathbf{\Sigma}}/\overline{\sigma}_{\mathbf{Q}}}t'}\right)$.\hfill$\blacksquare$

\subsection{Proof of Corollary \ref{corollary:centralized:stabilityNopReset}}
The arguments are similar to those used in the proof of Theorem \ref{theorem:centralized} by using the fact that in Lemma \eqref{proof:centralized:proposition:jumps} the expression in \eqref{eq:deltaVjumps} yields $\Delta V(\tilde{y}_c)\le 0$ whenever $\eta = 0$.

\vspace{-0.175cm}
\section{Conclusion}
\label{sec_conclusions}
In this paper, we explored decentralized concurrent learning dynamics with momentum and coordinated resetting for multi-agent systems over directed graphs. The proposed approach utilizes intermittent coordinated resets to enable collective convergence to a common parameter estimate, even with asymmetric information flow. Using Lyapunov theory for hybrid systems, we established input-to-state stability properties for the momentum-based dynamics, subject to a cooperative richness condition on the data matrices and a topology-dependent lower bound on the resetting frequency. We demonstrated the effectiveness of the proposed dynamics in cooperative adaptive control, showcasing their advantages in accelerated convergence and enhanced transient behavior compared to first-order adaptation algorithms.

\vspace{-0.2cm}
 \section*{Acknowledgments}
The authors would like to thank Bob Bitmead for fruitful discussions on transient performance and fundamental limitations in adaptive control and batch-based adaptive dynamics.

\vspace{-0.175cm}
\bibliographystyle{ieeetr}
\bibliography{references}

\begin{thebibliography}{10}

\bibitem{kamalapurkar2017concurrent}
R.~Kamalapurkar, B.~Reish, G.~Chowdhary, and W.~E. Dixon, ``Concurrent learning for parameter estimation using dynamic state-derivative estimators,'' {\em IEEE Trans. on Automatic Control}, vol.~62, no.~7, pp.~3594--3601, 2017.

\bibitem{vamvoudakis2015asymptotically}
K.~G. Vamvoudakis, M.~F. Miranda, and J.~P. Hespanha, ``Asymptotically stable adaptive--optimal control algorithm with saturating actuators and relaxed persistence of excitation,'' {\em IEEE transactions on neural networks and learning systems}, vol.~27, no.~11, pp.~2386--2398, 2015.

\bibitem{ochoa2021accelerated}
D.~E. Ochoa, J.~I. Poveda, A.~Subbaraman, G.~S. Schmidt, and F.~R. Pour-Safaei, ``Accelerated concurrent learning algorithms via data-driven hybrid dynamics and nonsmooth {ODE}s,'' {\em Learning for Dynamics and Control}, pp.~866--878, 2021.

\bibitem{casas2022switched}
J.~Casas, C.-H. Chang, and V.~H. Duenas, ``Switched adaptive concurrent learning control using a stance foot model for gait rehabilitation using a hybrid exoskeleton,'' {\em IFAC-PapersOnLine}, vol.~55, no.~41, pp.~187--192, 2022.

\bibitem{casas2023switched}
J.~Casas, C.-H. Chang, and V.~H. Duenas, ``Switched concurrent learning adaptive control for treadmill walking using a lower limb hybrid exoskeleton,'' {\em IEEE Trans. on Ctrl. Systems Technology}, 2023.

\bibitem{chowdhary2011theory}
G.~V. Chowdhary and E.~N. Johnson, ``Theory and flight-test validation of a concurrent-learning adaptive controller,'' {\em Journal of Guidance, Control, and Dynamics}, vol.~34, no.~2, pp.~592--607, 2011.

\bibitem{poveda2021data}
J.~I. Poveda, M.~Benosman, and K.~G. Vamvoudakis, ``Data-enabled extremum seeking: {A} cooperative concurrent learning-based approach,'' {\em International Journal of Adaptive Control and Signal Processing}, vol.~35, no.~7, pp.~1256--1284, 2021.

\bibitem{ochoa2022acceleratedADP}
D.~E. Ochoa and J.~I. Poveda, ``Accelerated continuous-time approximate dynamic programming via data-assisted hybrid control,'' {\em IFAC-PapersOnLine}, vol.~55, no.~12, pp.~561--566, 2022.

\bibitem{chowdhary2010concurrent}
G.~Chowdhary and E.~Johnson, ``Concurrent learning for convergence in adaptive control without persistency of excitation,'' in {\em 49th IEEE Conf. on Decision and Control (CDC)}, pp.~3674--3679, IEEE, 2010.

\bibitem{le2022concurrent}
J.~H. Le and A.~R. Teel, ``Concurrent learning in high-order tuners for parameter identification,'' in {\em 2022 IEEE 61st Conference on Decision and Control (CDC)}, pp.~2159--2164, IEEE, 2022.

\bibitem{nguyen2020momentumrnn}
T.~Nguyen, R.~Baraniuk, A.~Bertozzi, S.~Osher, and B.~Wang, ``Momentumrnn: Integrating momentum into recurrent neural networks,'' {\em Advances in Neural Information Processing Systems}, vol.~33, pp.~1924--1936, 2020.

\bibitem{muehlebach21}
M.~Muehlebach and M.~I. Jordan, ``Optimization with momentum: Dynamical, control-theoretic, and symplectic perspectives,'' {\em J. Mach. Learn. Res.}, vol.~22, Jan 2021.

\bibitem{su2014differential}
W.~Su, S.~Boyd, and E.~J. Cand\`{e}s, ``A differential equation for modeling nesterov's accelerated gradient method: Theory and insights,'' {\em J. Mach. Learn. Res.}, vol.~17, p.~5312–5354, jan 2016.

\bibitem{wibisono2016variational}
A.~Wibisono, A.~C. Wilson, and M.~I. Jordan, ``A variational perspective on accelerated methods in optimization,'' {\em proceedings of the National Academy of Sciences}, vol.~113, no.~47, pp.~E7351--E7358, 2016.

\bibitem{nguyen2022improving}
H.~H.~N. Nguyen, T.~Nguyen, H.~Vo, S.~Osher, and T.~Vo, ``Improving neural ordinary differential equations with nesterov's accelerated gradient method,'' {\em Advances in Neural Information Processing Systems}, vol.~35, pp.~7712--7726, 2022.

\bibitem{poveda2020heavy}
J.~I. Poveda and A.~R. Teel, ``The heavy-ball {ODE} with time-varying damping: Persistence of excitation and uniform asymptotic stability,'' in {\em 2020 American Control Conference (ACC)}, pp.~773--778, IEEE, 2020.

\bibitem{o2015adaptive}
B.~O’donoghue and E.~Cand\`{e}s, ``Adaptive restart for accelerated gradient schemes,'' {\em Foundations of computational mathematics}, vol.~15, no.~3, pp.~715--732, 2015.

\bibitem{roulet2017sharpness}
V.~Roulet and A.~d'Aspremont, ``Sharpness, restart and acceleration,'' {\em Advances in Neural Information Processing Systems}, vol.~30, 2017.

\bibitem{poveda2021robust}
J.~I. Poveda and N.~Li, ``Robust hybrid zero-order optimization algorithms with acceleration via averaging in time,'' {\em Automatica}, vol.~123, p.~109361, 2021.

\bibitem{wang2022scheduled}
B.~Wang, T.~Nguyen, T.~Sun, A.~L. Bertozzi, R.~G. Baraniuk, and S.~J. Osher, ``Scheduled restart momentum for accelerated stochastic gradient descent,'' {\em SIAM Journal on Imaging Sciences}, vol.~15, no.~2, pp.~738--761, 2022.

\bibitem{yu2020mass}
Y.~Yu, B.~A{\c{c}}{\i}kme{\c{s}}e, and M.~Mesbahi, ``Mass--spring--damper networks for distributed optimization in non-{E}uclidean spaces,'' {\em Automatica}, vol.~112, p.~108703, 2020.

\bibitem{boffi2020continuous}
N.~M. Boffi and J.-J.~E. Slotine, ``A continuous-time analysis of distributed stochastic gradient,'' {\em Neural computation}, vol.~32, no.~1, pp.~36--96, 2020.

\bibitem{sun2020continuous}
C.~Sun and G.~Hu, ``A continuous-time {Nesterov} accelerated gradient method for centralized and distributed online convex optimization,'' {\em arXiv preprint arXiv:2009.12545}, 2020.

\bibitem{ochoa2021momentum}
D.~E. Ochoa and J.~I. Poveda, ``Momentum-based {N}ash set-seeking over networks via multi-time scale hybrid dynamic inclusions,'' {\em IEEE Transactions on Automatic Control}, 2023.

\bibitem{ochoa2020robust}
D.~E. Ochoa, J.~I. Poveda, C.~A. Uribe, and N.~Quijano, ``Robust optimization over networks using distributed restarting of accelerated dynamics,'' {\em IEEE Ctrl. Systems Letters}, vol.~5, no.~1, pp.~301--306, 2020.

\bibitem{khong2014multi}
S.~Z. Khong, Y.~Tan, C.~Manzie, and D.~Ne{\v{s}}i{\'c}, ``Multi-agent source seeking via discrete-time extremum seeking control,'' {\em Automatica}, vol.~50, no.~9, pp.~2312--2320, 2014.

\bibitem{chen2006smooth}
X.~Chen and Y.~Li, ``Smooth formation navigation of multiple mobile robots for avoiding moving obstacles,'' {\em International Journal of Control, Automation, and Systems}, vol.~4, no.~4, pp.~466--479, 2006.

\bibitem{chen2013distributed}
W.~Chen, C.~Wen, S.~Hua, and C.~Sun, ``Distributed cooperative adaptive identification and control for a group of continuous-time systems with a cooperative pe condition via consensus,'' {\em IEEE Trans. on Automatic Control}, vol.~59, no.~1, pp.~91--106, 2013.

\bibitem{goebel2012hybrid}
R.~Goebel, R.~G. Sanfelice, and A.~R. Teel, {\em Hybrid dynamical systems}.
\newblock Princeton {U}niversity {P}ress, 2012.

\bibitem{sontag1995characterizations}
E.~D. Sontag and Y.~Wang, ``On characterizations of the input-to-state stability property,'' {\em Systems \& Control Letters}, vol.~24, no.~5, pp.~351--359, 1995.

\bibitem{bullo2018lectures}
F.~Bullo, {\em Lectures on network systems}, vol.~1.
\newblock CreateSpace, 2018.

\bibitem{cai2009characterizations}
C.~Cai and A.~R. Teel, ``Characterizations of input-to-state stability for hybrid systems,'' {\em Systems \& Ctrl. Ltrs.}, vol.~58, no.~1, pp.~47--53, 2009.

\bibitem{poveda2019codes}
J.~I. Poveda, K.~G. Vamvoudakis, and M.~Benosman, ``{CODES}: Cooperative data-enabled extremum seeking for multi-agent systems,'' in {\em 2019 IEEE 58th Conference on Decision and Control (CDC)}, pp.~2988--2993, IEEE, 2019.

\bibitem{javed2021excitation}
M.~U. Javed, J.~I. Poveda, and X.~Chen, ``Excitation conditions for uniform exponential stability of the cooperative gradient algorithm over weakly connected digraphs,'' {\em IEEE Control Systems Letters}, 2021.

\bibitem{wilson2021lyapunov}
A.~C. Wilson, B.~Recht, and M.~I. Jordan, ``A {L}yapunov analysis of accelerated methods in optimization.,'' {\em J. Mach. Learn. Res.}, vol.~22, pp.~113--1, 2021.

\bibitem{javed2021scalable}
M.~U. Javed, J.~I. Poveda, and X.~Chen, ``Scalable resetting algorithms for synchronization of pulse-coupled oscillators over rooted directed graphs,'' {\em Automatica}, vol.~132, p.~109807, 2021.

\bibitem{rios2017time}
H.~R{\'\i}os, D.~Efimov, J.~A. Moreno, W.~Perruquetti, and J.~G. Rueda-Escobedo, ``Time-varying parameter identiffication algorithms: {F}inite and fixed-time convergence,'' {\em IEEE {T}ransactions on Automatic Control}, vol.~62, no.~7, pp.~3671--3678, 2017.

\bibitem{tatari2021fixed}
F.~Tatari, M.~Mazouchi, and H.~Modares, ``Fixed-time system identification using concurrent learning,'' {\em IEEE {T}ransactions on Neural Networks and Learning Systems}, vol.~34, no.~8, pp.~1--11, 2021.

\bibitem{zang1990transient}
Z.~Zang and R.~R. Bitmead, ``Transient bounds for adaptive control systems,'' in {\em 29th IEEE conference on decision and control}, pp.~2724--2729, IEEE, 1990.

\bibitem{Khalil:1173048}
H.~K. Khalil, {\em {Nonlinear systems; 3rd ed.}}
\newblock Upper Saddle River, NJ: Prentice-Hall, 2002.

\bibitem{bookHDS}
R.~Goebel, R.~G. Sanfelice, and A.~R. Teel, {\em {Hybrid Dynamical Systems: Modeling, Stability and Robustness}}.
\newblock Princeton University Press, 2012.

\bibitem{gao2003exponential}
X.-B. Gao, ``Exponential stability of globally projected dynamic systems,'' {\em IEEE Transactions on Neural Networks}, vol.~14, no.~2, pp.~426--431, 2003.

\bibitem{zhang2015constructing}
H.~Zhang, Z.~Li, Z.~Qu, and F.~L. Lewis, ``On constructing {L}yapunov functions for multi-agent systems,'' {\em Automatica}, vol.~58, pp.~39--42, 2015.

\bibitem{bernstein2009matrix}
D.~S. Bernstein, {\em Matrix mathematics: theory, facts, and formulas}.
\newblock Princeton university press, 2009.

\bibitem{horn2012matrix}
R.~A. Horn and C.~R. Johnson, {\em Matrix analysis}.
\newblock Cambridge University Press, 2012.

\end{thebibliography}
\appendix
\label{sec:appendix_sec1}

\section{Auxiliary Lemmas}
\begin{lemma}\label{lemma:llt}
Consider the following block triangular matrix:
\begin{equation*}
    	M\coloneqq  \begin{pmatrix} A&B\\ 0&D\end{pmatrix}
\end{equation*}
  Suppose that $M$ is non-singular. Then, the minimum singular value of $M$, $\sigma_{\min}(M)$,
		satisfies
  \begin{equation*}
  \sigma_{\min}(M) \ge \frac{1}{\sqrt{\|A^{-1}\|^2(1+\|BD^{-1}\|^2)+\|D^{-1}\|^2}}    
  \end{equation*}
  
\end{lemma}
\textbf{Proof.} 
		First, since the inverse of the block triangular matrix $M$ is given by
		$$M^{-1}=\begin{bmatrix}A^{-1}&-A^{-1}BD^{-1}\\0&D^{-1} \end{bmatrix},$$
		we can upper-bound the 2-norm matrix of $M^{-1}$: 
        {\small
        \begin{align}
        \|M^{-1}\|^2&=\max _{|u|^2+|v|^2=1} \left |\begin{bmatrix}A^{-1}&-A^{-1}BD^{-1}\\0&D^{-1} \end{bmatrix}\begin{bmatrix}u\\v \end{bmatrix} \right |^2\notag\\
		&=\max _{|u|^2+|v|^2=1} \left |\begin{bmatrix}A^{-1}u-A^{-1}BD^{-1}v\\D^{-1}v \end{bmatrix} \right|^2\notag\\
		&=\max _{|u|^2+|v|^2=1} \left | A^{-1}u-A^{-1}BD^{-1}v \right| ^2+ \left |D^{-1}v \right|^2 \notag\\
        %
		&\le \|A^{-1}\|^2 (1+ \|BD^{-1}\|^2) + \|D^{-1}\|^2.	\label{aux:lemma1}
        \end{align}}
Then, since the minimum singular value of a matrix is the inverse of the 2-norm of the inverse matrix, i.e., $\sigma_{\min}(M)=\frac{1}{\|M^{-1}\|}$, we can use \eqref{aux:lemma1} to obtain the result.~\strut\hfill $\blacksquare$

\begin{lemma}\label{lem:Vbound}
For each $\tau_c\in[T_0,T]$ and $s\in\mathbb{R}_{\geq0}$, consider the following block matrix
\begin{align*}    
		\mathbf{V}_w(\tau_c,s)\coloneqq \begin{pmatrix}
  \frac{1}{\tau^2}\mathbf{Q}  &  \hat{\mathbf{\Omega}}(s)\\\vspace{0.2cm}  
  \hat{\mathbf{\Omega}}(s) &  \hat{\mathbf{\Sigma}}(s) \end{pmatrix},
\end{align*}
where 
\begin{align}
\hat{\mathbf{\Sigma}}(s)&\coloneqq(1-\omega)\mathbf{\Sigma}+k_{\text{t}} \mathbf{Q}  \mathbf{A}(s)\\
\hat{\mathbf{\Omega}}(s)&\coloneqq \mathbf{\Omega} + k_t \mathbf{Q}\mathbf{A}(s)
\end{align}
where $w\in[0,\omega]$, $\omega\in(0,1)$, and the matrices $\mathbf{Q} $, $\mathbf{\Omega}$, and $\mathbf{\Sigma}$ are as defined in \eqref{eq:sigma_omega}. Then, under Assumption \ref{assumption:CSR}, we have that:
\begin{equation}\label{eq:Xi_bound}
    \mathbf{V}_w(\tau_c,s) \succeq \nu I_{Nn},
\end{equation}
for all $\tau_c\in[T_0,T]$, all $w\in[0,\omega]$, and all $s\in\mathbb{R}_{\ge 0}$, where 
\begin{equation}
\nu:=\frac{(1-\omega)\underline{\sigma}_{\mathbf{\Sigma}}\underline{\sigma}_{\mathbf{Q}} - T^2   (\overline{\sigma}_{\mathbf{\Omega}}^{2}+k_t\chi^2)}{T^2((1-\omega)\underline{\sigma}_{\mathbf{\Sigma}}) + \underline{\sigma}_{\mathbf{Q}}}>0,
\end{equation}
with $\underline{\sigma}_{\mathbf{Q} }$, $\underline{\sigma}_{\mathbf{\Sigma}}$, and $\overline{\sigma}_{\mathbf{\Omega}}^{2}$ as defined in Proposition \ref{lemma:definitionQ}. \QEDB
\end{lemma}

\vspace{0.1cm}
\textbf{Proof:} First, we show that the matrix-valued function $\mathbf{V}_w(\cdot,\cdot)$ is positive-definite uniformly over $\tau_c\in[T_0,T]$, $s\in\mathbb{R}_{\ge}$, and $w\in[0,\omega]$. To do this, we decompose $\mathbf{V}_w$ as follows:
\begin{align}\label{eq:Xi_decompose}
\mathbf{V}_w(\tau_c,s)&= \mathbf{U}(\tau_c,s)^\top\mathbf{D}(\tau_c,s)\mathbf{U}(\tau_c,s)
\end{align}
where 
\begin{align*}
    \mathbf{D}(\tau_c,s)\coloneqq\begin{pmatrix} \frac{\mathbf{Q}}{\tau_c^{2}} &  0  \\ 0 & \hat{\mathbf{\Sigma}}(s)-\tau_c^{2}\hat{\mathbf{\Omega}}(s)^{\top}\mathbf{Q}^{-1}\hat{\mathbf{\Omega}}(s) \end{pmatrix},
\end{align*}
and 
\begin{align*}
\mathbf{U}(\tau_c,s)\coloneqq \begin{pmatrix}I &  \tau^{2}\mathbf{Q}^{-1}\hat{\mathbf{\Omega}}(s)\\ 0&I\end{pmatrix},
\end{align*}
Using the definition of $\mathbf{Q}$, and the fact that $\tau_c\leq T$ for all $\tilde{y}_c\in \mathbf{C}_c\cup \mathbf{D}_c$, we obtain
\begin{equation}\label{eq:Xi_diag1}
    \frac{\mathbf{Q} }{\tau^2}\succeq \left(\frac{\underline{\sigma}_{\mathbf{Q}}}{T^2}\right)I_{Nn}.
\end{equation}  
Also, it follows that
\begin{align}\label{eq:Xi_diag2}
		&\hat{\mathbf{\Sigma}}(s)- \tau_c^2 \hat{\mathbf{\Omega}}(s)^\top  \mathbf{Q} ^{-1}  \hat{\mathbf{\Omega}}(s)\succeq \zeta I_{Nn},
\end{align}
for all $s\in\mathbb{R}_{\ge0}$, where
\begin{equation}
\zeta:=(1-\omega)\underline{\sigma}_{\mathbf{\Sigma}} - \frac{ T^2}{\underline{\sigma}_{\mathbf{Q} }}(\overline{\sigma}_{\mathbf{\Omega}}^{2} + \chi^2(k_t)).
\end{equation}
Note that $\zeta>0$ since condition \eqref{theorem:centralized:frequencyband} holds by assumption. Therefore, since
\begin{align*}
     \frac{\mathbf{Q}}{\tau_c^2} \succ 0 ~~ \text{and} ~~ \hat{\mathbf{\Sigma}}(s) {-} \tau_c^2 \hat{\mathbf{\Omega}}^\top(s)  \mathbf{Q} ^{-1}  \hat{\mathbf{\Omega}}(s)\succ 0.
\end{align*}
it follows that  the matrix  $\mathbf{V}_w(\tau_c,s)$ is positive definite uniformly over $\tau_c\in[T_0,T]$, $s\in\mathbb{R}_{\ge 0}$, and $w\in[0,\omega]$ \cite[Theorem 7.7.7]{horn2012matrix}.

Now, we establish the matrix inequality \eqref{eq:Xi_bound}. To do so, we use the bounds \eqref{eq:Xi_diag1} and \eqref{eq:Xi_diag2} for \eqref{eq:Xi_decompose} to obtain that
\begin{align}
		\mathbf{V}_w(\tau_c,s)&\succeq \mathbf{U}^\top(\tau_c,s) 
  \begin{bmatrix}
  \frac{\underline{\sigma}_{\mathbf{Q} }}{T^2}I_{Nn}&0\\
  0&\zeta I_{Nn}	
  \end{bmatrix} \mathbf{U}(\tau_c,s)\notag\\
  &=\mathbf{V}(\tau_c,s)^\top \mathbf{V}(\tau_c,s),\label{eq:lowerXiV}
\end{align}
where $\mathbf{V}(\tau_c,s)$ is the upper block triangular matrix
\begin{align*}
\mathbf{V}(\tau_c,s)\coloneqq \begin{pmatrix}\sqrt{\frac{\underline{\sigma}_{\mathbf{Q}}}{T^{2}}}I_{Nn} & \sqrt{\frac{\tau^{4}\underline{\sigma}_{\mathbf{Q}}}{T^{2}}}\mathbf{Q}^{-1}\hat{\mathbf{\Omega}}(s)  \\  0 & \sqrt{\zeta} I_{Nn}\end{pmatrix}.
\end{align*}
By  applying Lemma \ref{lemma:llt} on the matrix $\mathbf{V}(\tau_c,s)$, and using \eqref{eq:lowerXiV} together with the fact that $\mathbf{V}$ has full column rank and thus that $\sigma_{\min}(\mathbf{V}^\top \mathbf{V})\ge\sigma_{\min}(\mathbf{V}^\top)\sigma_{\min}(\mathbf{V}) = \sigma^2_{\min}(\mathbf{V})$, we obtain
\begin{align*}
\mathbf{V}_w(\tau_c,s)& \succeq \frac{1}{\frac{T^2}{\underline{\sigma}_{\mathbf{Q}}}\left(1+\frac{\tau^4\underline{\sigma}_{\mathbf{Q}}}{\zeta T^2} \|\mathbf{Q}^{-1}\hat{\mathbf{\Omega}}(s) \|^2 \right) +\frac{1}{\zeta}} I_{2Nn} \\
  %
  %
  &= \frac{\zeta \underline{\sigma}_{\mathbf{Q}}^2}{T^2(\zeta \underline{\sigma}_{\mathbf{Q}} + T^2( \overline{\sigma}_{\mathbf{\Omega}}^2+\chi^2(k_t))) + \underline{\sigma}_{\mathbf{Q}}^2} I_{2Nn}\\
  %
  %
  &=\frac{(1-\omega)\underline{\sigma}_{\mathbf{\Sigma}}\underline{\sigma}_{\mathbf{Q}} - T^2   (\overline{\sigma}_{\mathbf{\Omega}}^{2}+\chi^2(k_t))}{T^2((1-\omega)\underline{\sigma}_{\mathbf{\Sigma}}) + \underline{\sigma}_{\mathbf{Q}}} I_{2Nn}
  %
  %
  %
\end{align*}	
where we have used the fact that the induced 2-norm is sub-multiplicative and that $\|\mathbf{Q}^{-1}\|\le 1/\underline{\sigma}_{\mathbf{Q}}$ and $\|\mathbf{\Omega}\|^2\le \overline{\sigma}_{\mathbf{\Omega}}^{2}$. This completes the proof. \hfill $\blacksquare$

\vspace{0.1cm}
\begin{lemma}\label{lemma:aux:ResettingPolicies}
Let $\eta:=(\eta_1,\eta_2,\dots,\eta_N)$ with $\eta_i\in\{0,1\}$ for all $i\in\mathcal{V}=\{1,2,\dots,N\}$ and $\mathbf{R}_\eta=\text{diag}(\eta)\otimes I_n$. Then, for all $\tilde{\theta},\tilde{p}\in\mathbb{R}^{Nn}$ we have: 
{\small
\begin{align*}
&|\mathbf{R}_\eta\tilde{\theta}+\left(I_{Nn}-\mathbf{R}_\eta\right)\tilde{p}-\tilde{\theta}|_{\mathbf{Q} }^2 + |\mathbf{R}_\eta\tilde{\theta}+\left(I_{Nn}-\mathbf{R}_\eta\right)\tilde{p}|_{\mathbf{Q} }^2 \\
&~~-|\tilde{p}|_{\mathbf{Q}}^2-|\tilde{p}-\tilde{\theta}|_{\mathbf{Q}}^2=|\tilde{\theta}|^2_{\mathbf{R}_\eta\mathbf{Q}} - |\tilde{p}|^2_{\mathbf{R}_\eta\mathbf{Q}} -  |\tilde{\theta} -\tilde{p}|^2_{\mathbf{R}_\eta\mathbf{Q}}
\end{align*}}
where $\mathbf{Q} $ is defined in \eqref{eq:Q_matrix}.\QEDB
\end{lemma}
\textbf{Proof:} By direct computation, we have:
\begin{align*}
 \big|\mathbf{R}_\eta\tilde{\theta}+\left(I_{Nn}-\mathbf{R}_\eta\right)\tilde{p}-\tilde{\theta}\big|_{\mathbf{Q} }^2   &=\big|\mathbf{R}_\eta(\tilde{\theta}-\tilde{p})-(\tilde{\theta}-\tilde{p})\big|_{\mathbf{Q} }^2\\
 &=|\left(I_{Nn}-\mathbf{R}_\eta\right)(\tilde{\theta}-\tilde{p})|^2_{\mathbf{Q}}\\
 &= |\mathbf{R}_{\eta}^c(\tilde{\theta}-\tilde{p})|^2_{\mathbf{Q}}\\
 &=|z|^2_{\mathbf{Q}},
\end{align*}
where $z\coloneqq\mathbf{R}_{\eta}^c(\tilde{\theta}-\tilde{p})$, and $\mathbf{R}_{\eta}^c\coloneqq I_{Nn}-\mathbf{R}_{\eta}$. . Writing $z=(z_1,\ldots, z_N)$, with $z_i=(\eta_i-1)\left(\tilde{\theta}_i-\tilde{p}_i\right)\in \mathbb{R}^n,~\forall i\in \mathcal{V}$, it follows that
\begin{align}
|z|^2_{\mathbf{Q}}
&=\sum_{i=1}^N q_i|\tilde{\theta}_i-\tilde{p}_i|^2(\eta_i-1)^2\notag\\&= \sum_{i=1}^N q_i|\tilde{\theta}_i-\tilde{p}_i|^2\left(1-\eta_{i}\right).
\label{lemma:reset1}
\end{align}
Similarly, 
\begin{align*}
|\mathbf{R}_\eta\tilde{\theta}+\left(I_{Nn}-\mathbf{R}_\eta\right)\tilde{p}|_{\mathbf{Q} }^2&=|\mathbf{R}_\eta(\tilde{\theta}-\tilde{p})+\tilde{p}|_{\mathbf{Q} }^2=|\tilde{z}|^2_{\mathbf{Q}},
\end{align*}
where $\tilde{z}=\mathbf{R}_\eta(\tilde{\theta}-\tilde{p})+\tilde{p}$. Writing $\tilde{z}\coloneqq(\tilde{z}_1,\ldots, \tilde{z}_N)$, with $\tilde{z}_i=\eta_i\left(\tilde{\theta}_i-\tilde{p}_i\right)+\tilde{p}_i\in \mathbb{R}^n,~\forall i\in \mathcal{V}$, we get:
\begin{align}
|\tilde{z}|_{\mathbf{Q}}^2
&=\sum_{i=1}^N q_i|\tilde{z}_i|^2\notag\\
&=\sum_{i=1}^Nq_i|\eta_i(\tilde{\theta}_i-\tilde{p}_i)+\tilde{p}_i|^2\notag\\
&=\sum_{i=1}^Nq_i\left(\eta_i^2|\tilde{\theta}_i-\tilde{p}_i|^2+2\eta_i(\tilde{\theta}_i-\tilde{p}_i)^\top(\tilde{p}_i)+|\tilde{p}_i|^2\right)\notag\\
%
%
&=\sum_{i=1}^N q_i\left(\eta_i(\tilde{\theta}_i-\tilde{p}_i)^\top(\tilde{\theta}_i+\tilde{p}_i) +|\tilde{p}_i|^2\right)\notag\\
&=\sum_{i=1}^Nq_i\eta_i|\tilde{\theta}_i|^2+\sum_{i=1}^N q_i|\tilde{p}_i|^2(1-\eta_i).\label{lemma:reset2}
\end{align}
Together \eqref{lemma:reset1} and \eqref{lemma:reset2} yield:
\begin{align*}
    |\mathbf{R}_\eta\tilde{\theta}+\big(I_{Nn}{-}\mathbf{R}_\eta&\big)\tilde{p}-\tilde{\theta}|_{\mathbf{Q} }^2 + |\mathbf{R}_\eta\tilde{\theta}+\left(I_{Nn}{-}\mathbf{R}_\eta\right)\tilde{p}|_{\mathbf{Q} }^2 \\ -|\tilde{p}|_{\mathbf{Q}}^2-|\tilde{p}-\tilde{\theta}|_{\mathbf{Q}}^2 
&= \sum_{i=1}^N q_i\eta_i |\tilde{\theta}_i|^2\\
&~~+ \sum_{i=1}^N q_i(1-\eta_i)\left(|\tilde{p}_i|^2 + |\tilde{\theta}_i -\tilde{p}_i|^2\right)\\
    &\quad-\sum_{i=1}^Nq_i|\tilde{p}_i|^2 -\sum_{i=1}^Nq_i|\tilde{\theta}_i -\tilde{p}_i|^2\\
    &~= \sum_{i=1}^N q_i\eta_i |\tilde{\theta}_i|^2\\
    &\qquad - \sum_{i=1}^Nq_i\eta_i\left(|\tilde{p}_i|^2 + |\tilde{\theta}_i -\tilde{p}_i|^2\right)\\
    &~= |\tilde{\theta}|^2_{\mathbf{R}_\eta\mathbf{Q}} - |\tilde{p}|^2_{\mathbf{R}_\eta\mathbf{Q}} -  |\tilde{\theta} -\tilde{p}|^2_{\mathbf{R}_\eta\mathbf{Q}}.
\end{align*}
\hfill $\blacksquare$

\end{document}